\documentclass[12pt]{article} 
\usepackage[sectionbib]{natbib}
\usepackage{array,epsfig,fancyheadings,rotating}
\usepackage[]{hyperref}  
\usepackage{sectsty, secdot}
\sectionfont{\fontsize{12}{14pt plus.8pt minus .6pt}\selectfont}
\renewcommand{\theequation}{\thesection\arabic{equation}}
\subsectionfont{\fontsize{12}{14pt plus.8pt minus .6pt}\selectfont}

\usepackage{booktabs}
\usepackage{amsmath}
\usepackage{amssymb}
\usepackage{amsfonts}
\usepackage{multirow}
\usepackage{amsthm}
\usepackage{xcolor}
\usepackage{mathrsfs}

\newcommand{\instructions}[1]{}

\usepackage{bm}

\newtheorem{theorem}{Theorem}

\newtheorem{corollary}{Corollary}
\newtheorem{proposition}{Proposition}
\theoremstyle{definition}
\newtheorem{definition}{Definition}

\newtheorem{remark}{Remark}
\usepackage{bm}
\usepackage{amsthm}

\defcitealias{FDA2006}{US FDA, 2006}
\defcitealias{FDA2018}{US FDA, 2018}

\def\pr{\textsf{P}} 
\def\ep{\textsf{E}} 
\def\Var{\mathrm{Var}}

\begin{document}


\renewcommand{\baselinestretch}{2}

\markright{ \hbox{\footnotesize\rm Statistica Sinica
}\hfill\\[-13pt]
\hbox{\footnotesize\rm
}\hfill }

\markboth{\hfill{\footnotesize\rm   LI-XIN ZHANG} \hfill}
{\hfill {\footnotesize\rm Efficient CAR} \hfill}

\renewcommand{\thefootnote}{}
$\ $\par


\fontsize{12}{14pt plus.8pt minus .6pt}\selectfont \vspace{0.8pc}
\centerline{\Large\bf Efficient Covariate-Adaptive
Randomization  }
\vspace{2pt} \centerline{\Large\bf  with Smallest Selection Bias in Clinical Trials}
\vspace{.4cm} \centerline{Li-Xin Zhang} \vspace{.4cm} \centerline{\it Zhejiang Gongshang University}   \vspace{.55cm} \fontsize{9}{11.5pt plus.8pt minus
.6pt}\selectfont


\begin{quotation}
\noindent {\it Abstract:} 
The power of conventional hypothesis tests and the selection bias in covariate-adaptive randomized clinical trials are typically investigated through simulation studies. In this article, we develop a theoretical framework for analyzing both the asymptotic power of linear-model-based tests for treatment effects and the asymptotic selection bias. Our results reveal that under covariate-adaptive randomization: (i) hypothesis tests generally lose power when covariates are not sufficiently balanced; in particular, the more covariates included in the testing model that are not accounted for in the randomization procedure, the greater the loss of power; (ii) the hypothesis test is usually more powerful than that under complete randomization; and (iii) compared with complete randomization, many popular covariate-adaptive randomization procedures in the literature—such as Pocock and Simon's marginal procedure, the stratified permuted block design, and Taves's minimization method—are generally efficient in terms of power but yield non-negligible selection bias. To address this trade-off, we propose a new family of covariate-adaptive randomization procedures that simultaneously account for both power and selection bias. Under these procedures, covariate imbalances are kept sufficiently small to achieve asymptotically maximal power for testing treatment effects, while the selection bias remains asymptotically negligible. These theoretical results provide a comprehensive picture of how the power of hypothesis testing, covariate imbalance, and selection bias interact with one another.

\vspace{9pt}
\noindent {\it Key words and phrases:}
Balancing covariates, clinical trial, loss of power, selection bias, Pocock and Simon's procedure
\par
\end{quotation}\par

\def\thefigure{\arabic{figure}}
\def\thetable{\arabic{table}}

\renewcommand{\theequation}{\thesection.\arabic{equation}}

\fontsize{12}{14pt plus.8pt minus .6pt}\selectfont

\setcounter{section}{1} 
\setcounter{equation}{0} 

\lhead[\footnotesize\thepage\fancyplain{}\leftmark]{}\rhead[]{\fancyplain{}\rightmark\footnotesize\thepage}
\setcounter{section}{0}
\section{Introduction}
\label{s:intro}
\setcounter{equation}{0}

It is well known that covariates often play an important role in clinical trials. Clinical trialists frequently express concern over unbalanced treatment arms with respect to key covariates of interest. To balance important covariates, covariate-adaptive designs \citep{RosenbergerLachin2002} are usually employed. The most commonly implemented methods are stratified randomization and marginal minimization \citep{McEntegart2003}. Stratified randomization  first stratifies the population and then applies separate randomization within each stratum;  examples include the stratified permuted block design, among others. To handle many covariates, \cite{Taves1974} and \cite{PocockSimon1975} introduced the marginal minimization method, which attempts to minimize the weighted sum of  differences of allocation numbers between  treatment groups across all marginal values of covariates. \cite{HuHu2012} discussed some limitations of these classical designs, proposed a generalized family of covariate-adaptive designs,  and established  their theoretical properties. For   further discussion on handling covariates in clinical trials, see \cite{McEntegart2003}, \cite{RosenbergerSverdlov2008}, \cite{HuRosenberger2006}, and  references therein.

It is important to study the statistical inference associated with covariate-adaptive designs. In practice, conventional tests are often employed without accounting for the covariate-adaptive randomization scheme. It remains a concern whether conventional tests are still valid under covariate-adaptive designs, especially when the covariates used in trial design and those incorporated in inference procedures are not the same. A substantial body of simulation literature has examined statistical inference under covariate-adaptive designs. \cite{Birkett1985} and \cite{Forsythe1987} raised concerns about  the validity of unadjusted analysis under covariate-adaptive designs and suggested, based on simulation studies, that all covariates used in Taves's minimization method should be included in the analysis to achieve a valid test. \cite{FeinsteinLandis1976} and \cite{GreenByar1978} studied inference problems for stratified randomization with binary responses. \cite{Ciolino2011} showed that power loss could be  non-negligible if the balancing  of important continuous covariates is  ignored  even when  adjustment is made in the analysis  for those covariates. More discussions can be found in \cite{Simon1979}, \cite{TuShalayPater2000}, \cite{Aickin2009} and elsewhere. \cite{ShaoYuZhong2010} and \cite{MaHuZhang2015} theoretically proved that, when some or all covariates used in trial design are not incorporated in inference procedures, the hypothesis  test for comparing treatment effects is usually conservative in terms of Type I error. It is now generally accepted that covariates used in trial design should also be incorporated in inference procedures.

In this paper, we consider the opposite  direction, namely, that covariates incorporated in inference procedures should be used in trial design. We aim to study the exact power of the hypothesis test for comparing treatment effects.  An expansion of the relative power function is established. It is  shown that when some or all covariates used in inference procedures are not incorporated in trial design, the  test generally loses power; the more covariates in the testing model   that are not incorporated in the randomization procedure, the greater the power loss. However,  adaptive randomization is usually more powerful than complete randomization. When all the covariates incorporated in inference procedures are   sufficiently balanced, the power of the  test  is asymptotically equivalent to the  maximum achievable power.  A broad class of covariate-adaptive designs, which includes most of the covariate-adaptive designs in the literature, for example, Pocock and Simon's marginal procedure, the stratified permuted block design,  can achieve the maximum power.

In the literature, continuous covariates are typically discretized in order to be included in the randomization scheme \citep{Taves2010}. However, as discussed in \cite{Scott2002},  breaking  a continuous covariate into subcategories entails increased effort and loss of information. \cite{Ciolino2011} also pointed out that  a lack of practical methods for balancing continuous covariates, combined with limited knowledge of  the cost of failing to balance them, has led to a common phenomenon whereby continuous covariates are excluded from the randomization plan in clinical trials. Our theoretical results on  power illustrate how  discretization affects  power loss and provide a formula for the cost of failing to balance continuous covariates.

On the other hand, randomization in clinical trials   is fundamental to study design. Randomization is desirable for a number of reasons including the mitigation of  selection bias which may arise if the person in charge of selecting patients  has advance knowledge of the treatment assignments.  \cite{Efron1971} discussed how to measure the lack of randomness and treatment imbalance, and proposed a biased coin design (BCD) to   trade off between the treatment imbalance and lack of randomness. However, the lack of randomness of covariate-adaptive randomization has been seldom considered theoretically. Based on the measure of selection bias defined by \cite{Efron1971}, we   find that, compared to complete randomization, the selection bias of both Pocock and Simon's marginal procedure and stratified permuted block design is  non-negligible, even though the asymptotic power of these designs can be maximal. A  natural question is whether there exists a covariate-adaptive randomization procedure under which both  power and  selection bias are simultaneously optimized. We propose a new family of covariate-adaptive designs that include Pocock and Simon's marginal procedure, \cite{HuHu2012}'s design, and \cite{Wei1978}'s adaptive biased coin design. Within this family, a new covariate-adaptive randomization procedure is defined by choosing a suitable allocation function, under which covariate imbalances are small enough that the power of testing treatment effects is asymptotically maximal, while at the same time the asymptotic selection bias attains the minimum value of $1/2$, so that every guessing strategy is asymptotically equally useless against the allocation procedure.

In Section \ref{s:testing}, we study the power of the hypothesis  tests for comparing the treatment effects. Section \ref{s:selectionbias}  addresses the theory of  selection bias. In Section \ref{s:procedure}, the new family of adaptive randomization procedures  are proposed. The proofs are provided in the Supplementary Material.

\section{The power of hypothesis test under covariate-adaptive designs}
\label{s:testing}
\setcounter{equation}{0}

In this section, we study the power of hypothesis testing based on a linear model framework for covariate-adaptive designs. Suppose two treatments $1$ and $2$ are studied under a covariate-adaptive randomized clinical trial, where $\mu_1$ and $\mu_2$ are parameters measuring the mean effects of treatment $1$ and $2$, respectively. Let $n$ be the total number of patients enrolled in the study. Let $T_i$ be the assignment of the $i$th patient, i.e., $T_i=1$ for treatment 1 and $T_i=0$ for treatment 2, $i=1,2,...,n$. Conditional on the treatment assignment $T_i$, the following linear model is assumed for the response of the $i$th patient $Y_i$,
\begin{equation}\label{model}
Y_i=\mu_1T_i+\mu_2(1-T_i)+\beta_1 X_{i,1} +...+ \beta_I X_{i,I} +\varepsilon_i,
\end{equation}
where
  \begin{enumerate}
  \item[-] $X_{i,k}$s  are discrete or continuous covariates which are independent and identically distributed as $X_k$ , $k=1,\ldots, I$;
  \item[-] some or all values of $X_{i,k}$s  are used in the randomization procedure, $k=1,..., I$;
  \item[-] all covariates are independent of each other, and $EX_k=0$ and $\Var(X_k)>0$ for all $k$, $k=1,\ldots, I$;
  \item[-] $\varepsilon_i$s are independent and identically distributed normal random errors with zero mean   and variance $\sigma_{\varepsilon}^2$ and are independent of $X_k$,$k=1,\ldots, I$.
\end{enumerate}
Note that  $X_{i,k}$   are assumed to be scalars in model \eqref{model}. If $X_{i,k}$   is a discrete covariate, $X_{i,k}$   is a scalar that can take several values corresponding to different categories. In practice, a vector is usually used to represent a discrete covariate with multiple categories.  In this paper, $X_{i,k}$  is assumed to be a scalar for simplicity, but all the results can be extended to the situation where discrete covariates with multiple categories are represented by vectors.

We let $\bm Y=(Y_1, Y_2, ..., Y_n)^T$, $\bm{\varepsilon}=(\varepsilon_1, \varepsilon_2, ..., \varepsilon_n)^T$, $\bm\beta=(\beta_1,..., \beta_I)^T$, $\bm\gamma=(\mu_1, \mu_2, \beta_1,..., \beta_I)^T$, and
$$\bm X=\left[\begin{array}{ccccc} T_1 & 1-T_1 & X_{1,1}&\cdots&X_{1,I}\\
                                  T_2 &1-T_2 &X_{2,1}&\cdots&X_{2,I}\\
                                  \vdots& \vdots&\vdots&\ddots&\vdots\\
                                 T_n &1-T_n &X_{n,1}&\cdots&X_{n,I}\end{array}\right].
$$
Then model \eqref{model} can be written as
$
\bm Y=\bm X\bm\gamma+\bm{\varepsilon}.
$
The ordinary least squares method is used to obtain the estimate of $\bm\gamma$, which has the explicit form
\begin{align*}
\hat{\bm \gamma}=(\bm X^T\bm X)^{-1}\bm X^T\bm Y
=\bm \gamma+(\bm X^T\bm X)^{-1}\bm X^T\bm \varepsilon.
\end{align*}
In the covariate-adaptive randomization, the results of allocations $T_1,\ldots, T_n$ depend only on the covariates $\bm X_1,\ldots, \bm X_n$, and so are independent of $\bm\epsilon$. Hence, given $\bm X$, $\hat{\bm\gamma}$ follows a multi-dimensional normal distribution with mean $\bm \gamma$ and variance-covariance matrix $\sigma_{\epsilon}^2(\bm X^T\bm X)^{-1}$.

When model \eqref{model} is constructed to study data from a covariate-adaptive randomized clinical trial, the primary interest is usually to compare treatment effects between different groups.  To compare treatment effects of $\mu_1$ and $\mu_2$, the following right-sided hypothesis test is usually used
\begin{equation}\label{test}
H_0: \mu=0  \mbox{ versus } H_A: \mu > 0.
\end{equation}
where $\mu=\mu_1-\mu_2$ is the treatment effect difference.  The
right-sided $t$-test    statistic for \eqref{test} is
\begin{equation}\label{stat}
T=\frac{\bm L \boldsymbol{\hat\gamma}}{(\widehat\sigma_{\epsilon}^2\bm L(\boldsymbol{X}^T\boldsymbol{X})^{-1}\bm L^T)^{1/2}},
\end{equation}
where    $\bm L=(1, -1, 0,..., 0)$, and $\hat\sigma_{\epsilon}^2=(\bm Y-\bm X\hat{\bm \gamma})^T(\bm Y-\bm X\hat{\bm \gamma})/(n-I-2)$ is the estimate of $\sigma^2_{\epsilon}$.  If $T>t_{1-\alpha}(\nu)$, where $t_{1-\alpha}(\nu)$ is the  $({1-\alpha})$th quantile of a $t$-distribution with the degree of freedom $\nu=n-I-2$, we  reject the null hypothesis, otherwise we accept the null hypothesis. Given $\bm X$, the conditional power function is
$$\beta_{T,n}(\mu|\bm X)=1-F\left(t_{1-\alpha}(\nu); \nu, \frac{\mu}{\sigma_{\epsilon}^2\bm L(\bm X^T\bm X)^{-1}\bm L^T)^{1/2}}\right), $$
 where $F(t;\nu,\delta)$ is the distribution function of a non-central $t$-distribution with non-central parameter $\delta$ and degree of freedom  $\nu =n-I-2$.

 Let $\bm X_i=(X_{i,1},\cdots, X_{i,I})^T$,
\begin{align*}
    N_{n,1}=\sum_{i=1}^n T_i, \;        N_{n,2}= & \sum_{i=1}^n (1-T_i), \;   D_n=N_{n,1}-N_{n,2},  \\
     \bm D_n^x=\sum_{i=1}^n (2T_i-1)\bm X_i, \;      D_n^{x_k}& = \sum_{i=1}^n (2T_i-1)X_{i,k},\; k=1,\ldots, I, \\
    \overline{\bm X}_{n,1}=\frac{\sum_{i=1}^nT_i \bm X_i}{N_{n,1}},\; &
    \; \overline{\bm X}_{n,2}=\frac{\sum_{i=1}^n(1-T_i) \bm X_i}{N_{n,2}}, \\
    \overline{\bm Y}_{n,1}=\frac{\sum_{i=1}^nT_i \bm Y_i}{N_{n,1}},\; &
    \; \overline{\bm Y}_{n,2}=\frac{\sum_{i=1}^n(1-T_i) \bm Y_i}{N_{n,2}},
\end{align*}
\begin{align*}
    \bm S_{n,xx}=& \sum_{i=1}^n T_i(\bm X_i-\overline{\bm X}_{n,1})(\bm X_i-\overline{\bm X}_{n,1})^T \\
    &+\sum_{i=1}^n (1-T_i)(\bm X_i-\overline{\bm X}_{n,2})(\bm X_i-\overline{\bm X}_{n,2})^T.
    \end{align*}
Here $D_n$ is the   difference between the numbers of patients in treatment group 1 and 2 overall, which can be regarded as a measure of overall treatment imbalance, and,
$\bm D_n^x=\sum_{i=1}^n T_i\bm X_i-\sum_{i=1}^n (1-T_i)\bm X_i$ can be regarded as a measure of covariate imbalance. It can be seen that
$$\overline{\bm X}_{n,1}=\frac{n\overline{\bm X}_n+\bm D_n^x}{n+D_n}\; \text{ and } \;  \overline{\bm X}_{n,2}=\frac{n\overline{\bm X}_n-\bm D_n^x}{n-D_n} $$
are functions of $D_n$ and $\bm D_n^x$, where $\overline{\bm X}_n=\sum_{i=1}^n \bm X_i/n$. Also
\begin{align*}
   \bm L(\bm X^T\bm X)^{-1}\bm L^T
 =& \frac{1}{N_{n,1}}+\frac{1}{N_{n,2}}+(\overline{\bm X}_{n,1}-\overline{\bm X}_{n,2})^T\bm S_{xx}^{-1}(\overline{\bm X}_{n,1}-\overline{\bm X}_{n,2}) \\
 =& \frac{4}{n\big(1-(D_n/n)^2\big)}+\frac{4Q_n^2}{n\big(1-(D_n/n)^2\big)^2},
 \end{align*}
where
$$ Q_n^2= \big(\bm D_n^x/\sqrt{n}-\overline{\bm X}_nD_n/\sqrt{n}\big)^T\left(\frac{\bm S_{n,xx}}{n}\right)^{-1}\big(\bm D_n^x/\sqrt{n}-\overline{\bm X}_nD_n/\sqrt{n}\big). $$
Hence, the    power function can be written as
$$\beta_{T,n}(\mu|\bm X)=1-F\left(t_{1-\alpha}(\nu); \nu, \frac{\mu}{2\sigma_{\epsilon}}\ell_n\right), $$
where
$$ \ell_n= \sqrt{\frac{n\big(1-(D_n/n)^2\big)^2}{1- (D_n/n)^2+\frac{1}{n}Q_n}}\le \sqrt{n}. $$
Because $F(t;
 \nu,\delta)$ is a non-increasing function of $\delta$, it is obvious that when $D_n=0$ and $\bm D_n^x=\bm 0$,  the power function takes its largest value
 $$\beta_{T,n}(\mu|0)=1-F\left(t_{1-\alpha}(\nu); \nu, \frac{\mu}{2\sigma_{\epsilon}}\sqrt{n}\right). $$
A clinical trial  is only a single realization of a random phenomenon, and it cannot be assumed that the observed imbalances $D_n$ and $\bm D_n^x$ will be zero. Therefore, the loss of power is unavoidable. The problem is, compared to the largest power, how small is the realized power acceptable?  The distance between the power functions $\beta_{T,n}(\mu|\bm X)$ and $\beta_{T,n}(\mu|0)$ may reflect  the loss of power. However,
since the power functions $\beta_{T,n}(\mu|\bm X)$ and $\beta_{T,n}(\mu|0)$ are both convergent  to $1$, from their difference one  can not discern how the imbalances $D_n$ and $\bm D_n^x$ affect the loss of power clearly. To deal with this problem,   we consider the relative loss of power as the ratio:
$$ LossP_{T,n}(\mu|\bm X)=\frac{1-\beta_{T,n}(\mu|0)}{1-\beta_{T,n}(\mu|\bm X)} =\frac{F\left(t_{1-\alpha}(\nu); \nu, \frac{u}{2\sigma_{\epsilon}}\sqrt{n}\right)}{F\left(t_{1-\alpha}(\nu); \nu, \frac{\mu}{2\sigma_{\epsilon}}\ell_n\right)},  $$
$ \nu=n-I-2$. The smaller   value of $LossP_{T,n}(\mu|\bm X)$, the greater  the power loss.
The following theorem gives the order of the loss of power.
\begin{theorem}\label{thPower} Write $\sigma_{x,k}^2=\Var(X_k)$, $k=1,\ldots, I$,
$$V_n= (D_n/\sqrt{n})^2+\sum_{k=1}^I \frac{(D_n^{x_k}/\sqrt{n})^2}{\sigma_{x,k}^2}. $$
 Suppose $\frac{D_n}{n}$ and $\frac{\bm D_n^x}{n}$ converge to zero  in probability (or almost surely), then
\begin{equation}\label{eqLopT}LossP_{T,n}(\mu|\bm X)= \exp\left\{ -\frac{\mu^2}{8\sigma_{\epsilon}^2}V_n\big(1+o(1)\big)\right\}.
\end{equation}
in probability (or almost surely).
\end{theorem}

In Theorem \ref{thPower}, the value of $V_n=(D_n/\sqrt{n})^2+\sum_{i=1}^p(D_n^{x_k}/\sqrt{n})^2\sigma_{x,k}^{-2}$ gives the order of the loss of power. For  complete randomization in which each patient is allocated to treatment $1$ or $2$ with the same probability $1/2$,
$$ \left(\frac{D_n}{\sqrt{n}},\frac{D_n^{x_1}}{\sqrt{n}},\cdots,\frac{D_n^{x_I}}{\sqrt{n}}\right)^T\overset{d}\to N_{I+1}\left(\bm 0, diag(1,\sigma_{x,1}^2,\cdots,\sigma_{x,I}^2)\right), $$
 where $N_{I+1}$ is an $(I+1)$-dimensional normal distribution with mean and variance-covariance matrix as specified. So, $V_n\overset{d}\to \chi^2(I+1)$, where $\chi^2(I+1)$ is a $\chi^2$-distributed random variable with  $I+1$  degree of freedom. Hence we have
 \begin{corollary}\label{CorCR} If the patients are allocated to treatments by  the complete randomization, then
 \begin{equation}\label{eqPowerCR.2} LossP_{T,n}(\mu|\bm X)\overset{d}\to \exp\left\{-\frac{\mu^2}{8\sigma^2}\chi^2(I+1)\right\}<1.
 \end{equation}
\end{corollary}
By \eqref{eqPowerCR.2}, the expected relative loss of power is
$$  \ep\left[LossP_{T,n}(\mu|\bm X)\right]\cong   \left(1+\frac{\mu^2}{4\sigma_{\epsilon}^2}\right)^{-\frac{I+1}{2}}.  $$

 The right-side of \eqref{eqLopT} will converge to 1 whenever $V_n\to 0$. Therefore, we conclude the following corollary.
  \begin{corollary}\label{Cor2}
 If
\begin{equation} \label{eqthPowerC1.1}\frac{D_n}{\sqrt{n}}\to 0 \; \text{  and }\; \frac{\bm D_n^x}{\sqrt{n}}\to 0
\end{equation}
 in probability (or almost surely), then
 $  LossP_{T,n}(\mu|\bm X) \to 1 $
in probability (or almost surely).
\end{corollary}
\begin{definition}
When the condition  \eqref{eqthPowerC1.1} is satisfied,   the power of the hypothesis test for comparing treatment effects   under this allocation procedure is equivalent to that under the completely balanced allocation which has the largest power. Such a randomization procedure is called {\em an efficient covariate-adaptive design}.
\end{definition}
Condition \eqref{eqthPowerC1.1} means that the covariates used in the testing model should be sufficiently balanced. In clinical trials, covariate-adaptive designs are usually based on discrete covariates (Taves 2010).
If a continuous (or general) covariate is to be used in randomization, a discrete conversion need to be performed to break down the continuous covariate into a discrete variable with several subcategories. In general,  the covariate-adaptive design is applied with respect to discrete variables $d_k(X_k)$, where $d_k$s are discrete functions. In such cases, define $\delta_k= X_k-E[X_k|d_k(X_k)]$ and $\tilde X_{k}=d_k(X_k)$.
Accordingly,  $\delta_{i,k}= X_{i,k}-E[X_{i,k}|d_k(X_{i,k})]$ and $\tilde X_{i,k}=d_k(X_{i,k})$ are  $i$th observations of covariates $\delta_k$ and $\tilde X_k$, $k=1,\ldots, I$.    Suppose $ \tilde X_{k}$ have $m_k$ levels, resulting in $\prod_{k=1}^{I}m_k$ strata in total. Let $\tilde{\bm X}_i=(\tilde X_{i,1},\ldots,\tilde X_{i,I})$ represent the covariate profile of the $i$th patient used in randomization, i.e., $\tilde{\bm X}_i=(\tilde{x}_1^{t_1}, \tilde{x}_2^{t_2},\ldots,\tilde{x}_I^{t_I})$ if ${\tilde X_{i,k}}$ is at level $\tilde{x}_k^{t_k}$. For convenience, we use $(t_1, \ldots,t_I)$ to denote the stratum formed by patients who have the same covariate profile $(\tilde{x}_1^{t_1},\ldots,\tilde{x}_I^{t_I})$, and use $(j;t_j)$ to denote the margin formed by patients with $ \tilde X_{j} =\tilde{x}_j^{t_j}$. Then let
\begin{description}
  \item[-] $D_n=N_{n,1}-N_{n,2}=\sum_{i=1}^n (2T_i-1)$ be the difference between the numbers of patients in treatment group 1 and 2 as total;
  \item[-] $D_n(j;t_j)=\sum_{i=1}^n(2T_i-1)I\{\tilde{X}_{i,j}=\tilde{x}_j^{t_j}\}$  be the differences between the numbers of patients in the two treatment groups on the margin $(j;t_j)$;
  \item[-]  $D_n(t_1,\ldots,t_I)=\sum_{i=1}^n (2T_i-1)I\{\tilde{\bm X}_i=(\tilde{x}_1^{t_1},\ldots,\tilde{x}_I^{t_I})\}$ be the differences between the numbers of patients in the two treatment groups  within the stratum $(t_1,\ldots,t_I)$.
\end{description}
These differences play important roles in the covariate-adaptive randomization procedure  and in properties of statistical inference for covariate-adaptive designs. Each one of these differences is used to measure the \emph{imbalance} at the corresponding level (overall, marginal, or within-stratum) of the  $n$ assignments.

\begin{theorem}\label{thPower2} Let $\sigma_{\delta,k}^2=\ep[Var(\delta_k|d_k(X_k))]=\ep[\delta_k^2]$, $k=1,\ldots,p$. Assume that each assignment $T_m$ depends only on $T_1,\ldots, T_{m-1}$ and $\tilde{\bm X}_1,\ldots,\tilde{\bm X}_m$, i.e., the probability of $T_m=1$ is a function of $T_1,\ldots, T_{m-1}$ and $\tilde{\bm X}_1,\ldots,\tilde{\bm X}_m$.  Suppose
\begin{equation}\label{eqthPowerC2.1}\frac{D_n}{\sqrt{n}}\to 0\; \text{ and } \frac{D_n(j;t_j)}{\sqrt{n}}\to 0
\end{equation}
  in probability (or almost surely), $t_j=1,\ldots, m_j$, $j=1,\ldots,I$. Then
\begin{equation}\label{eqthPower2.T}
LossP_{T,n}(\mu|\bm X)
\overset{d}\to \exp\left\{-\frac{\mu^2}{8\sigma^2_{\epsilon}}\sum_{k=1}^I \frac{\sigma_{\delta,k}^2}{\sigma_{x,k}^2}\chi_k^2(1)\right\},
\end{equation}
where $\chi_1^2(1),\ldots,\chi_I^2(1)$ are independent $\chi^2$-distributed   random variables  with $1$ degree of freedom. The expected relative loss of power is
$$  \ep\left[LossP_{T,n}(\mu|\bm X)\right]\cong   \prod_{i=1}^I\left(1+\frac{\mu^2}{4\sigma_{\epsilon}^2}\frac{\sigma_{\delta,k}^2}{\sigma_{x,k}^2}\right)^{-1/2}. $$
\end{theorem}

 Condition \eqref{eqthPowerC2.1} is   satisfied under various
covariate-adaptive designs for balancing covariates. For example, under stratified permuted block designs,
  the class of covariate-adaptive designs proposed by \cite{HuHu2012}, and   \cite{PocockSimon1975}'s marginal
procedure,
\begin{equation} \label{b:probability}D_n \; \text { and } \; D_n(j;t_j),\; \text{ are bounded in probability},
\end{equation}
$t_j=1,\ldots, m_j$, $j=1,\ldots, I$ (c.f. \cite{MaHuZhang2015,HuYeZhang2023}), and so \eqref{eqthPowerC2.1} is satisfied.  Further, if the covariates $X_k$s are discrete and  $d_k(X_k)=X_k$, $k=1,\ldots, I$, then $D_n$ and $\bm D_n^{x}$ are bounded in probability, resulting in  $V_n=O(n^{-1})$ and
$$ V_n=O(n^{-1}) \; \text{ and }  LossP_{T,n}(\mu|\bm Z)=1+O(n^{-1}) \text{ in probability}.$$
 Notice that $O(n^{-1})$ is the fastest convergence rate unless the treatments are completely balanced so that $D_n=0$ and $D_n(j,t_j)=0$ for all $t_j$ and $j$, which indicates that Efron's biased coin design, as well as its generalizations such as   Pocock and Simon's marginal procedure and Hu and Hu's design, is usually efficient in terms of power with respect to the problem it addresses. 

 Simulation studies have shown that \cite{Taves1974}'s minimization method also reduces marginal imbalances as well as the overall imbalance. However, theoretical evidence has not, to the best of our knowledge, been found in literature. Whether \eqref{b:probability} holds or not is still an open problem. We will show that \eqref{eqthPowerC2.1} is true under Taves's minimization method (c.f., Theorem \ref{theorem1} and Corollary \ref{Cor3}), so it is an  efficient covariate-adaptive design in the sense that $ LossP_{T,n}(\mu|\bm X) \overset{P}\to 1$  if the covariates $X_k$s are discrete and  $d_k(X_k)=X_k$, $k=1,\ldots, I$.

In Theorem \ref{thPower2}, the ratio $\sigma_{\delta,k}^2/ \sigma_{x,k}^2$ describes  the cost to   loss of power of failing to balance the $k$th covariate completely.  It is obvious that $0\le \sigma_{\delta,k}^2/ \sigma_{x,k}^2\le 1$. The larger  $\sigma_{\delta,k}^2$ is, the larger the cost is. If $\ep[X_k|d_k(X_k)]=X_k$, i.e., the $k$-covariate is not considered in the randomization,  then $\sigma_{\delta,k}^2/ \sigma_{x,k}^2=1$ and accordingly, the $k$-covariate  contributes fully to the loss of power and so the cost is the largest. If $\ep[X_k|d_k(X_k)]=d_k(X_k)$, i.e., the $k$th covariate is discrete and is balanced enough  with respect to all of its values, then $\sigma_{\delta,k}^2=0$ and accordingly, the contribution of the $k$-covariate to the loss of power can be neglected. So, the fewer  covariates  and their values  that are balanced, the greater the power loss. Comparing the results in Corollary \ref{CorCR} and Theorem \ref{thPower2}, we find that the hypothesis test for comparing treatment effects under a covariate-adaptive randomization procedure is usually more
powerful than the one under the complete randomization.

In the literature, covariate-adaptive randomization schemes  are usually based on discrete covariates. When the covariates are continuous,  the covariates are usually discretized.   However, as discussed in \cite{Scott2002}, the discretization  means increased effort and loss of information. \cite{Ciolino2011} also pointed out that uncontrolled imbalances in continuous covariates may affect the statistical analysis of the trial outcome. In the formula \eqref{eqthPower2.T}, $\sigma_{\delta, k}^2$ is the variance of the uncontrolled part of the $k$th covariate. The ratio  $\sigma_{\delta,k}^2/ \sigma_{x,k}^2$ just gives the cost  of failing to balance the  covariate $X_k$ completely. It is easily seen that the cost can be decreased by refining
the discretization.  Recently, \cite{Ma2024} proposed a procedure to balance the continuous covariates $\bm X_i$s directly  so  that the condition \eqref{eqthPowerC1.1} is satisfied under widely satisfied conditions. Under \cite{Ma2024}'s procedure, the asymptotic power is the largest when the covariates considered in the randomization and the covariates used in the hypothesis test
for treatment effects are consistent. 

The inconsistency of the covariates considered in the randomization and the covariates used in the hypothesis test
for treatment effects may affect the hypothesis test. The evidences have been found in many simulation studies  \citep{Birkett1985, Forsythe1987, TuShalayPater2000,Aickin2009,Ciolino2011}.  \cite{MaHuZhang2015}, \cite{Zhu2025} etc derived the theory for  the case   where the covariates considered in the testing are fewer than those considered in the randomization under a large class of covariate-adaptive randomization procedures satisfying  \eqref{eqthPowerC2.1}.   Combining the results in Theorem \ref{thPower2} and those of  \cite{MaHuZhang2015}, we conclude that, under a large class of covariate-adaptive designs for balancing covariates, (i) if the covariates considered in the testing are less than those considered in the randomization, then the hypothesis test  is usually conservative in terms of  Type I error, and, the more of the values of covariates are not consistent, the more the testing is conservative; (ii) if the covariates considered in the testing are more than those considered in the randomization, then the hypothesis test usually loses power, and, the more of the values of covariates are not consistent, the more   the power loss; (iii) in any case, the hypothesis test  is usually more
powerful than the one under  complete randomization.

\begin{remark} \label{remark2.1} The conclusions of Theorems \ref{thPower} and \ref{thPower2} remain true when the statistic $|T|$ is used to test the following two-sided hypothesis test
\begin{equation}\label{twosidetest}
H_0: \mu=0  \mbox{ versus } H_A: \mu \ne  0.
\end{equation}
When $\sigma_{\epsilon}^2$ is known, Theorems \ref{thPower} and \ref{thPower2} are also valid for one-sided or two-sided $z$-test.
\end{remark}

The formula \eqref{eqLopT} remains true under the
alternative hypothesis with $0<\mu=\mu_n\to 0$ and $\mu_n\sqrt{n}\to \infty$. But it no longer holds when a sequence of local
alternative hypotheses of the type
$H_A:\mu =\eta/\sqrt{n}>0$ are considered. In such case, one can show that $\frac{d}{d\delta}F(t;\nu,\delta)$ converges to a continuous function $f_{\infty}(t;\delta)<0$ as $\nu \to \infty$ uniformly in $t$ and $\delta$ on bounded intervals.  It follows that
\begin{align*}
&F\Big(t_{1-\alpha}(\nu);\nu,\frac{\eta}{2\sigma_{\epsilon}}\frac{\ell_n}{\sqrt{n}}\Big)-F\Big(t_{1-\alpha}(\nu);\nu,\frac{\eta}{2\sigma_{\epsilon}}\Big)\\
 =&-\frac{\eta}{4\sigma_{\epsilon}}f_{\infty}\Big(z_{1-\alpha};\frac{\eta}{2\sigma_{\epsilon}}\Big) \frac{n-\ell_n^2}{n} \big(1+o(1)\big)=C_{\alpha,\eta} \frac{V_n}{n}\big(1+o(1)\big).
 \end{align*}
 Hence
\begin{equation}\label{eq:power:local} \beta_{T,n}(u|0)-\beta_{T,n}(u|\bm X)=  O(1) \frac{V_n}{n}.
\end{equation}
The value of $V_n$ also reflects the order of  power loss.

 When the distribution of the model error $\epsilon_i$ is not normal, though it is difficult to derive the exact power function,
  \eqref{eq:power:local} is expected to be true especially when the errors follow a continuous distribution.  An approximation  method is usually applied to approximate  the true distribution of the statistic  $T$ for hypothesis testing, when the distribution of error $\epsilon_i$ is not normal or unknown.
 It is well known that if normal approximation is used, the approximation precision of the distribution as well as the power function is $O\left(1/\sqrt{n}\right)$.   If the Bootstrap method is used, the approximation precision  would be improved to $O\left(1/n\right)$. Hence, compared to the approximation precision, the loss of power $O(1) \frac{V_n}{n}$ can be ignored whenever $V_n\to 0$. Again, the condition \eqref{eqthPowerC1.1} is needed.

\section{Selection bias}
\label{s:selectionbias}
\setcounter{equation}{0}
  
 As advocated in \cite{Efron1971}, selection bias (lack of randomness)   in a design manifests in the possibility for the experimenter
to  partially guess the sequence of treatment allocations; thus it can be measured, say, by the
expected percentage of correct guesses when an optimal guessing strategy is used. Let $J_m =1$ if the $m$th assignment
is guessed correctly, and $J_m=0$ otherwise. The expected proportion of correct guesses is
$$SB_n(\Delta)=\ep\left[\frac{1}{n}\sum_{m=1}^n J_m\right]=\frac{1}{n}\sum_{m=1}^n\pr(J_m=1). $$
Here $SB_n(\Delta)$  stands for the selection bias of the allocation procedure $\Delta$.
Every guessing strategy
is equally useless against   complete randomization, yielding the expected percentage
of correct guesses $\frac{1}{2}$  for any $n$; this is   the optimal value of $SB_n$. 

For  
adaptive randomization, the optimal strategy consists, at each step,
of picking the under-represented treatment, which is always the treatment with the higher probability
of being allocated, with no preference in case of a tie. Let $p_m$ be the conditional probability of assigning the $m$th patient to treatment $1$ given the history $\sigma$-field including the information of past allocations, past covariates and current covariate observed. The probability that the $m$th assignment
is guessed correctly is $\ep[p_m\vee (1-p_m)]$. So,
$$SB_n(\Delta)=\frac{1}{n}\sum_{m=1}^n \ep\left[p_m\vee(1-p_m)\right]=\frac{1}{2}+ \frac{1}{n}\sum_{m=1}^n \ep\Big[\big|p_m-\frac{1}{2}\big|\Big].$$
There are other ways  to measure the selection bias or predictability of a randomization procedure. For example, an equivalent way  is to consider \cite{Smith1984}'s index, namely the difference
between the mean percentage of correct guesses and that of incorrect guesses:
$U_n(\Delta)=2SB_n(\Delta)-1,$
which has the advantage of ranging from 0 to 1. Also, the selection bias can be expressed
as  $\frac{1}{n}\sum_{m=1}^n\ep[|p_m-1/2|]$, 
which is just the \cite{Smith1984}'s index $U_n(\Delta)=2SB_n(\Delta)-1$.   
  For the discussion on the predictability and selection bias in     randomized block procedure, one can refer to \cite{BergerIvanovaKnoll2003}. In this paper, we only consider Efron's definition. 
It is easily seen that $SB_n(\Delta)>1/2$ unless $p_m\equiv 1/2$ for all $m$. In general, we consider the asymptotic selection bias $SB(\Delta)=\lim_{n\to \infty}SB_n(\Delta)$.

Under stratified permuted block designs, or \cite{Taves1974}'s minimization method, a positive percentage of patients are assigned deterministically, and so are guessed correctly, resulting in
$$ \liminf_{n\to \infty} SB_n>1/2. $$
For the class of biased coin designs, when the design does not include the covariate,\cite{Efron1971} gives the formula for the asymptotic selection bias $SB$ of his original biased coin design, which indicates that $SB(BCD)>1/2$.

For \cite{PocockSimon1975}'s  procedure, no closed form of the asymptotic selection bias $SB$ has been found in the literature. But we have the following result.

 \begin{theorem} For the Pocock and Simon's procedure, \cite{HuHu2012}'s procedure or more general procedures as in Remark \ref{remark4.1}, we have
 $$ \lim_{n\to \infty} SB_n>\frac{1}{2}. $$
 \end{theorem}

 So, for all the covariate-adaptive randomization procedures that have been introduced so far, when the condition   \eqref{eqthPowerC2.1} is satisfied so that the loss of power would be negligible, the asymptotic selection bias $SB$ is larger than $1/2$ so that the lack of randomness is non-negligible.

  It is obvious that, if
  \begin{equation}\label{eq:p:convergence} p_m\to \frac{1}{2}\;\; \text{ in probability},
  \end{equation}
  then $SB$ attains the minimum value $1/2$  of the selection biases, and every guessing strategy
is asymptotically  equally useless against the design.
  When the trials do not include the covariate, \cite{Wei1978} defines the adaptive  BCD satisfying \eqref{eq:p:convergence}  such that $SB=1/2$. But Wei's design does not satisfy \eqref{eqthPowerC2.1}. The problem is whether we can find a randomization procedure such that  $SB=1/2$ and the condition   \eqref{eqthPowerC2.1} is also satisfied? The next section will deal with this problem.

\section{A New family of Covariate-Adaptive Randomization methods}
\label{s:procedure}
\setcounter{equation}{0}

In this section, we consider the randomization methods.
We first give a general framework for covariate-adaptive randomization. We consider the same  setting  as that of \cite{PocockSimon1975} and only focus on two treatment groups $1$ and $2$ here.
As before, let $T_{j}$ be the assignment of the $j$th patient, $j=1,\ldots,n$,
i.e., $T_{j}=1$ for treatment 1 and $T_{j}=0$ for treatment 2.
Recall that $\tilde{\bm X}_j=\big(d_1(X_{j,1}),\ldots, d_I(X_{j,I})\big)$ indicates the discrete part of the  covariate profile of that patient considered in adaptive randomization,
i.e., $\tilde{\bm X}_j=(\tilde{x}_1^{t_1},\ldots,\tilde{x}_I^{t_I})$ if his or her $i$th covariate is at level $\tilde{x}_i^{t_i}$, $1\leq i \leq I $ and $1\leq t_i\leq m_i$.
For convenience, we use $(t_1,\ldots,t_I)$ to denote the \emph{stratum} formed by patients who possess the same covariate profile $(\tilde{x}_1^{t_1},...,\tilde{x}_I^{t_I})$, resulting in $M=\prod_{i=1}^{I}m_i$ strata, and use $(i;t_i)$ to denote the \emph{margin} formed by patients whose $i$th covariate is at level $\tilde{x}_i^{t_i}$. Recall that
$D_m=\sum_{j=1}^m (2T_j-1)=N_{m,1}-N_{m,2}$, $D_{m}(i;t_i)=\sum_{j=1}^m(2T_j-1)I\{\tilde{X}_{j,i}=\tilde{x}_i^{t_i}\}$  and 
$D_m(t_1,\ldots,t_I)=\sum_{j=1}^m (2T_j-1)I\{\tilde{\bm X}_j=(\tilde{x}_1^{t_1},\ldots,\tilde{x}_I^{t_I})\}$ are measures of the overall, marginal and within-stratum imbalances  of the first $m$ assignments, respectively.

\noindent{\bf 4.1 The randomization procedure}

The  procedure is defined as follows:
\begin{enumerate}
\item[1)] The first patient is assigned to treatment 1 with probability 1/2.

\item[2)] Suppose $(n-1)$ patients have been assigned to  treatments ($n>1$) and
the covariate value $\bm X_n=\tilde{\bm x}_n^{\ast}=(\tilde{x}_1^{t_1^{\ast}},\ldots,\tilde{x}_I^{t_I^{\ast}})$ of the $n$th patient is observed and falls within  stratum $(t_1^{*},\ldots,t_I^{*})$.

 For the first $(n-1)$ patients, let $D_{n-1}$, $D_{n-1}(i;t_i^{\ast})$ and $D_{n-1}(t_1^{\ast},\ldots,t_I^{\ast})$ be the overall imbalance, the marginal imbalance for  the margin $(i;t_i^{*})$ and the within-stratum imbalance for the stratum $(t_1^{\ast},\ldots,t_I^{\ast})$, respectively. 

\item[3)] If the $n$th patient were assigned to treatment $1$, then $D^{(1)}_{n}=D_{n-1}+1$ would
be the ``potential'' overall difference in the two groups; similarly,
\begin{displaymath}
D^{(1)}_{n}(i;t_i^{\ast})= D_{n-1}(i;t_i^{\ast})+1
\end{displaymath}
and
\begin{displaymath}
D^{(1)}_{n}(t_1^{\ast},\ldots,t_I^{\ast})=D_{n-1}(t_1^{\ast},\ldots,t_I^{\ast})+1
\end{displaymath}
would be the potential differences on  margin $(i;t_i^{\ast})$ and within  stratum  $(t_1^{\ast},\ldots,t_I^{\ast})$, respectively.

\item[4)] Define an imbalance measure $Imb_{n}^{(1)}$ by
\begin{displaymath}
Imb_{n}^{(1)}=w_{o}[D^{(1)}_n]^2+\sum_{i=1}^{I}w_{m,i}[D^{(1)}_n (i;t_i^{\ast})]^2
+w_{s}[D^{(1)}_n(t_1^{\ast},\ldots,t_I^{\ast})]^2,
\end{displaymath}
      which is the weighted imbalance that would be caused if the $n$th patient
      were assigned to treatment 1. $w_{o}$, $w_{m,i}$ $(i=1,\ldots,I)$ and $w_{s}$ are nonnegative
      weights placed on overall, within a covariate margin and within a stratum
      cell, respectively. Without loss of generality we can assume
\begin{displaymath}
w_{o}+w_{s}+\sum_{i=1}^{I}w_{m,i}=1.
\end{displaymath}
\item[5)] In the same manner we can define $Imb_{n}^{(2)}$, the weighted imbalance that would
be caused if the $n$th patient were assigned to treatment 2. In this case, the three types
of potential differences are the existing ones minus 1, instead of plus 1.

\item[6)]  Conditional on the assignments of the first $(n-1)$ patients as well as the covariates' profiles of the first $n$ patients, assign the $n$th patient to treatment 1 with probability
\begin{align}\label{eqallocationP}
&P(T_{n}=1|\bm{Z}_{n-1},\tilde{\bm X}_n=(\tilde{x}_1^{t_1^{\ast}},\ldots, \tilde{x}_I^{t_I^{\ast}}),\bm{T}_{n-1})\nonumber\\
& \quad = g_n\left(Imb_{n}^{(1)}-Imb_{n}^{(2)}\right)
\end{align}
  where $n>1$,  $\bm{Z}_{n-1}=(\tilde{\bm X}_1,\ldots,\tilde{\bm X}_{n-1})$,
  $\bm{T}_{n-1}=(T_1,\ldots,T_{n-1})$.  Here  $g_n(x)$ is a   real  function with $0\le g_n(x)\le 1$. It is called an allocation function.
  \end{enumerate}

Using the basic equation $(x+1)^2-(x-1)^2=4x$, the critical quantity $Imb_{n}^{(1)}-Imb_{n}^{(2)}$ in Step $6)$  can be simplified as
\begin{align}
&Imb_{n}^{(1)}-Imb_{n}^{(2)}\nonumber\\
=&4\left\{w_{o}D_{n-1}+\sum_{i=1}^{I}w_{m,i}D_{n-1}(i;t_i^{\ast})+w_{s}D_{n-1}(t_1^{\ast},\ldots,t_I^{\ast})\right\}\nonumber\\
:=&4\cdot\Lambda_{n-1}(t_1^{\ast},\ldots,t_I^{\ast}) \label{simplification}
\end{align}
Therefore, the allocation probability $g_n\left(Imb_{n}^{(1)}-Imb_{n}^{(2)}\right)$ is determined by the value of $\Lambda_{n-1}(t_1^{\ast},\ldots,t_I^{\ast})$, which is a weighted average of current imbalances at different levels.

This framework includes various kinds of covariate-adaptive randomization procedures. If  stratified randomization is considered, one  need only to let $w_s=1$ and the other weights  be zero. If $w_o=w_s=0$, the design is a marginal randomization procedure.  

\noindent{\bf 4.2 The choice of the allocation function}

  In the literature different views have been given as to the selection of the allocation probability function $g_n(\cdot)$. \cite{Efron1971}, \cite{PocockSimon1975}, \cite{HuHu2012} suggested
\begin{equation} \label{eqBCDfunction} g_n(x)=  1-\rho  \text{ if } x>0,\;
\frac{1}{2}    \text{ if } x=0,\;
\rho    \text{ if } x<0,
\end{equation}
where $1/2<\rho<1$. If the allocation function is chosen in this way,  and $w_s=w_o=0$, then the design is the \cite{PocockSimon1975}'s marginal procedure. \cite{MaHuZhang2015}  theoretically proved that under Pocock and Simon's   procedure, the marginal imbalances and overall imbalance  are bounded in probability. But the selection bias   is   not negligible  as we have seen in Section \ref{s:selectionbias}.

 When the design does not include the covariate, \cite{Wei1978} defines the adaptive  BCD with allocation function  defined as
 \begin{equation}\label{WeiBCD}
g_n(x)=g\left(\frac{x}{n-1}\right),
\end{equation}
 where $g(x)$ is  a non-increasing function with $0\le g(x)\le 1$, $g(x)=1-g(-x)$,  $g(0)=0.5$ and  $g^{\prime}(0)<0$. It has been shown that
 $ SB(WeiBCD)=1/2, $
which attains the minimum value   as   complete randomization does.
However,  Wei's BCD does not satisfy the condition \eqref{eqthPowerC2.1} and so the power loss could be non-negligible  \citep{AntogniniGiovagnoli2004, Antognini2008}.

Next, we show that  a   covariate-adaptive design can be defined  by choosing  suitable allocation function $g_n(x)$  so that the selection bias is asymptotically at its minimum  and the covariate imbalances considered are  of the order of $o(n^{1/2})$ in probability for which the loss of power would be negligible.

In general,   we assume   the allocation function $g_n(x)$ ($0\le g_n(x)\le 1$) satisfies the following conditions
\begin{gather}
    g_n(x)\le 1/2 \le g_n(-x)  \text{ when } x\ge 0, \label{eqCondF2}\\
    \frac{\left|1/2-g_n(x_n)\right|}{|x_n|/n}\to +\infty \; \text{ whenever }\; 0<\frac{|x_n|}{n}\to 0\label{eqCondF3}
   \end{gather}
  and
  \begin{equation}\label{eqCondF4}
  g_n(x_n)\to \frac{1}{2} \text{ whenever }\; \frac{x_n}{n}\to 0.
  \end{equation}
We will show that if the allocation function satisfies the conditions \eqref{eqCondF2} and \eqref{eqCondF3}, then covariate imbalances  are  of the order of $o(n^{1/2})$ in probability and,  if the conditions
 \eqref{eqCondF2} and \eqref{eqCondF4} are satisfied, then the asymptotic selection bias is $1/2$.
 The allocation function $g_n(x)$ of BCD given in \eqref{eqBCDfunction} satisfies \eqref{eqCondF2} and \eqref{eqCondF3}, but it does not satisfy condition \eqref{eqCondF4}. Wei's allocation function $g_n(x)$ given in \eqref{WeiBCD} satisfies    \eqref{eqCondF2} and \eqref{eqCondF4} , but it does not satisfy condition \eqref{eqCondF3}.

 It is easily seen that
$g_n(x)=1-\Phi\left(sgn(x)\sqrt{|x|/n}\right)$ satisfies all  conditions \eqref{eqCondF2}-\eqref{eqCondF4}, where $\Phi(x)$ is the standard normal distribution function. It can be chosen as an allocation function to define a design as desired.

In the following we will give a specific  allocation function. For this allocation function, the design has many fine properties.
Let $0\le g(x)\le 1$ be a non-increasing function with $g(0)=0.5$,  $g^{\prime}(0)<0$. Define
\begin{equation}\label{eqalloacFuct}
g_n(x)=g\left(\frac{x}{(n-1)^{\gamma}}\right), \;\; 0\le  \gamma\le 1.
\end{equation}
When $\gamma=1$, the allocation function $g_n(x)$ is Wei's function given in \eqref{WeiBCD}. When $\gamma=0$,  the design is a generalization of the  \cite{PocockSimon1975}'s procedure as well as  \cite{HuHu2012}'s design. In this paper, we mainly consider the case $0<\gamma<1$.

 \noindent{\bf 4.3 Theoretical Properties}

To consider the  asymptotic properties of the design, we need some more notations. For the first $n$ patients, we know that $D_n(t_1,\ldots,t_I)$
is the true difference between the two treatment arms   within stratum $(t_1,\ldots,t_I)$. Let
$$
\bm {D}_{n}=\left[D_n(t_1,\ldots,t_I)\right]_{1\leq t_1\leq m_1,\ldots, 1\leq t_I\leq m_I}
$$
be an array of dimension $m_1\times\ldots \times m_I$ which stores the current assignment
differences in all strata and so stores the current imbalances.
Also,  the covariates $\tilde{\bm X}_1,\tilde{\bm X}_{2},\ldots$ are independently and identically
distributed. Since $\tilde{\bm X}_n=(\tilde{x}_1^{t_1},...,\tilde{x}_I^{t_I})$ can take
 $M=\prod_{i=1}^{I} m_i$  different values, it  in fact follows an $M$-dimension
multinomial distribution with parameter $\bm p=(p(t_1,\ldots,t_I))$, where each element   $p(\bm t)=\pr(\tilde{\bm X}_n=(\tilde{x}_1^{t_1},\ldots, \tilde{x}_I^{t_I}))$ is the probability that a patient falls within the corresponding stratum.  Without loss of generality, we assume $p(t_1,\ldots,t_I)>0$ for all $(t_1,\ldots,t_I)$.

 Write $\bm t=(t_1,\ldots,t_I)$. Define
 $$ M_n=\sum_{\bm t}w_s D_n^2(\bm t)+\sum_{i=1}^I\sum_{t_i=1}^{m_i} w_{m,i} D_n^2(i;t_i)+w_o D_n^2. $$

Now we give our main results.
\begin{theorem}\label{theorem1}
 Suppose the allocation function $g_n(x)$ satisfies the conditions \eqref{eqCondF2} and \eqref{eqCondF3}.
Then
\begin{equation}\label{eqth1.1}
\ep\left(\frac{ M_n }{n}\right)^r\to 0 \;\;  \forall r>0,
\end{equation}
i.e., $M_n=o(n)$ in $L_r$ for all $r>0$.  In particular,
\begin{description}
  \item[\rm (i)] If $w_s>0$, then $\bm D_n=o(n^{1/2})$  in $L_r$  for any $r>0$;
  \item[\rm (ii)]  If $w_s+w_{m,i}>0$, then $D_n(i;t_i)=o(n^{1/2})$  in $L_r$  for any $r>0$, $t_i=1,\ldots, m_i$;
  \item[\rm (iii)]  For any case $D_n=o(n^{1/2})$  in $L_r$  for any $r>0$.
\end{description}
Further, if the condition \eqref{eqCondF4} is satisfied, then
$ SB_n \to \frac{1}{2}  $ 
and   $\bm D_n=o(n)$ in $L_r$  for any $r>0$.
\end{theorem}

\begin{theorem}\label{theorem2} Let the allocation function $g_n(x)$ be defined as \eqref{eqalloacFuct} with $0<\gamma\le 1$, and $0\le g(x)\le 1$ be a non-increasing function with $g(0)=1/2$,  $g^{\prime}(0)<0$.
Then
\begin{equation}\label{eqth2.1}
M_n=O(n^{\gamma}) \;\; \text{ in } L_r   \text{ for all } r>0,
\end{equation}
\begin{equation}\label{eqth2.2}
M_n=o(n^{\gamma+\epsilon}) \; a.s.    \text{ for any } \epsilon>0
\end{equation}
and
\begin{equation}\label{selectionBias}
SB_n =\frac{1}{2}+ O(n^{-\gamma/2}).
\end{equation}
 In particular,
\begin{description}
  \item[\rm (i)] If $w_s>0$, then $\bm D_n=O(n^{\gamma/2})$  in $L_r$  for any $r>0$, and $\bm D_n=o(n^{\gamma/2+\epsilon})$ a.s. for any $\epsilon>0$;
  \item[\rm (ii)]  If $w_s+w_{m,i}>0$, then $D_n(i;t_i)=O(n^{\gamma/2})$  in $L_r$  for any $r>0$, and $D_n(i;t_i)=o(n^{\gamma/2+\epsilon})$ a.s. for any $\epsilon>0$, $t_i=1,\ldots, m_i$;
  \item[\rm (iii)]  For any case $D_n=O(n^{\gamma/2})$ in $L_r$  for any $r>0$, and $D_n=o(n^{\gamma/2+\epsilon})$ a.s. for any $\epsilon>0$, $\bm D_n=o(n)$ a.s. and in   $L_r$  for any $r>0$;
  \item[\rm (iv)] Suppose $w_s+w_{m,i}>0$, $X_i=d_i(X_i)$ (the covariates in the test and in the randomization procedure are the same), $i=1,\ldots, I$. Then
  $$ LossP_{T,n}(\mu|\bm X)=1+O\big(n^{\gamma-1}\big) \;\;\text{ in probability}. $$
\end{description}
Further, in \eqref{selectionBias}, $O(n^{-\gamma/2})$ can be written as $ a_n n^{-\gamma/2}$ with
$$c_0<\liminf_{n\to \infty} a_n\le  \limsup_{n\to \infty}a_n\le  \frac{2\sqrt{ |g^{\prime}(0)|}}{ 2-\gamma}, $$
where $c_0>0$ is a constant which does not depend on $\gamma$.
\end{theorem} 

The value of $\gamma$ gives the order of the allocation bias as well as the order of selection bias, which provides a picture of how the allocation bias and the selection bias conflict with each other.  The smaller  $\gamma$ is, the smaller  the allocation bias is, but the larger the selection bias is.
In practice, one may choose $\gamma=1/2$ to give a tradeoff between the two types of bias, both of which are moderate.

\begin{remark}\label{remark4.1} Theorem \ref{theorem2} does not include the case of $\gamma=0$. When $\gamma=0$, suppose $0<g(x)\le 1/2\le g(-x)<1$ for all $x\ge 0$, and $\liminf\limits_{x\to \infty}|g(x)- 1/2|>0$. It can be shown that $\left[\Lambda_n(t_1,\ldots,t_I)\right]_{t_i=1,\ldots, m_i, i=1,\ldots,I}$ is a positive recurrent Markov chain with an invariant    distribution $\bm \pi$ and $\sup_n\ep\left[M_n^r\right]<\infty$ for all $r>0$ (c.f.  \cite{HuZhang2020}). Further,
$$ \lim_{n\to \infty}SB_n=\frac{1}{2}+\sum_{\bm t} \ep_{\bm \pi}\Big[\Big|g\big(\Lambda(\bm t)\big)-\frac{1}{2}\Big|\Big]p(\bm t)>1/2.$$
So, Theorem \ref{theorem2} is still valid for $\gamma=0$. 
\end{remark}

It should be noticed  that  \cite{Taves1974}'s minimization method is not included in Remark \ref{remark4.1}, because the related allocation function is $g_n(x)=1/2$ if $x=0$, $g_n(x)=0$ if $x>0$ and $g_n(x)=1$ if $x<0$.   This allocation function $g_n(x)$ satisfies the conditions \eqref{eqCondF2} and \eqref{eqCondF3}. So \eqref{eqth1.1} and (i)-(iii) of Theorem \ref{theorem1} remain true. We have  more precise results as in the following corollary.
\begin{corollary}\label{Cor3}
 Suppose the allocation function $g_n(x)$ satisfies the conditions \eqref{eqCondF2} and there exist  constants $x_0>0$ and $p_0>0$ such that
\begin{equation}\label{eqcor1.1} \big|g_n(x)-\frac{1}{2}\big|\ge p_0 \;\; \text{ for all } |x|\ge x_0 \text{ and } n\ge 1.
\end{equation}
Then
\begin{equation}\label{eqcor1.2}
M_n=O(n^{\epsilon}) \; a.s. \; \text{ and } \;\; \text{ in } L_r   \text{ for all } r>0, \epsilon>0.
\end{equation}
\end{corollary}
\section{Concluding Remarks}
\label{s:conclusion}
\setcounter{equation}{0}
 
In this paper, we have studied the theoretical properties  of  hypothesis testing and selection bias under  covariate-adaptive designs. First, based on a linear normal model framework, we derive the corresponding asymptotic power  of the hypothesis test for comparing treatment effects.   The asymptotic relative power is a decreasing function of the covariate imbalances, characterizing the impact of covariate imbalance on statistical power. Second, it is shown that, for many classical covariate-adaptive randomization procedures
 including  stratified permuted block designs,  \cite{Taves1974}'s minimization method, \cite{PocockSimon1975}'s  procedure, and  \cite{HuHu2012}'s procedure, the selection bias is  non-negligible. Finally,  a new family  of covariate-adaptive designs are proposed  under which the power of testing the treatment effects is  asymptotically the largest  and, at the same time,  the asymptotic selection bias attains the minimum value. 

To apply Theorem \ref{thPower2} on the asymptotic power to a specific covariate-adaptive randomized clinical trial and find an efficient  randomization, one  only needs  to check the condition \eqref{eqthPowerC2.1} or \eqref{eqthPowerC1.1}. The conclusion can be applied to a broad range of designs.  To check or find the most random covariate-adaptive randomization procedure, one  only needs  to check a simple  condition \eqref{eq:p:convergence}.  The new family of covariate-adaptive randomization procedures are shown to satisfy the condition \eqref{eqthPowerC2.1} and \eqref{eq:p:convergence} simultaneously. The results in this paper provide new insights about balance, selection bias,  efficiency, and randomization methods  in clinical trials, and the framework can be used to study properties of other statistical methods under covariate-adaptive randomization procedures and to define new covariate-adaptive randomization procedures.

The condition \eqref{eqthPowerC1.1} is important for ensuring that the asymptotic loss of  power can be neglected so that the randomization procedure is efficient. When the covariates are discrete, this condition is satisfied for a large family of covariate-adaptive designs, including many popular designs in the literature, and the new family proposed in this paper.  However, all
covariate-adaptive designs considered in this article are based
on discrete covariates. When the covariates are continuous,  they are typically discretized in order to be included in the randomization scheme.
 Formula \eqref{eqthPower2.T} gives  the cost in terms of power loss of failing to balance the $k$th covariate completely when a  randomization scheme is implemented by breaking down a continuous covariate into a discrete variable with several subcategories.   \cite{Ma2024} proposed a procedure that can balance the continuous covariates $\bm X_i$s directly. Let $D_n=\sum_{i=1}^n (2T_i-1)$ and $\bm D_n^x=\sum_{i=1}^n (2T_i-1)\bm X_i$ be defined as in Section \ref{s:testing}. Denote the measure of imbalance by $M_n=(D_n)^2+\|\bm D_n^x\|^2$, $Imb_n^{(1)}=M_n\big|_{T_n=1}$ and $Imb_n^{(2)}=M_n\big|_{T_n=0}$. Ma et al's procedure assigns  the $n$th patient to treatment 1 with probability $\rho$ ($>1/2$) when $Imb_n^{(1)}<Imb_n^{(2)}$,
$1/2$ when $Imb_n^{(1)}=Imb_n^{(2)}$, and $1-\rho$ when $Imb_n^{(1)}>Imb_n^{(2)}$. Ma et al. proved  that condition \eqref{eqthPowerC1.1} is satisfied under some suitable conditions.   This procedure can be generalized by defining the allocation probability as $g_n(Imb_n^{(1)}-Imb_n^{(2)})$. We may have similar results to those in Theorems \ref{theorem1}-\ref{theorem2} and Corollary \ref{Cor3}. In particular, if $g_n(x)$ satisfies the conditions \eqref{eqCondF2}-\eqref{eqCondF4} in Theorem \ref{theorem1}, we have
$$ LossP_{T,n}(\mu|\bm X) \to 1 \text{ in probability and } \lim_{n\to \infty} SB_n=\frac{1}{2}. $$

The theoretical properties for covariate-adaptive
designs can be generalized in several ways.
First, the theory  and adaptive randomization procedures are based on clinical trials with two treatments,
but they can be generalized to multiple treatments \citep{HuZhang2004,TymofyeyevRosenbergerHu2007, HuYeZhang2023}. Second, a key
assumption to derive theoretical results on the power is the independence
between covariates. We may apply a similar idea to dependent
covariates by incorporating the correlation structure.     Third, to derive the relative loss of power, another important assumption is the normality of the regression errors.
  It is an open problem whether the conclusions in Theorems \ref{thPower} and \ref{thPower2} are still valid for a general distribution of errors, especially when an approximate distribution is applied. This may not be an easy problem, but precise self-normalized Cram\'er-type
large deviations are helpful for solving it (c.f. \cite{JingShaoWang2003}). Fourth,   the randomization procedures considered  in this paper only depend on the covariates;   patient  responses to treatments are not incorporated. It is possible to generalize the results to efficient response-adaptive designs \citep{HuZhangHe2009} and covariate-adjusted adaptive designs \citep{Zhang2007}.
These topics are left for future
research.

\section*{Supplementary Material}

The online supplementary material contains the theoretical proof of all the theorems, and the simulation study of the new CAR with a specific allocation function $g$ in \eqref{eqalloacFuct}.

\section*{Acknowledgments}
The authors thank the Associate Editor and the two anonymous referees
for many helpful comments and suggestions. The research was partly   supported by  grants from  NSF of China (No. U23A2064), National Key R\&D Program of China (No. 2024YFA1013502)  and the Summit Advancement Disciplines of Zhejiang Province (Zhejiang Gongshang University - Statistics).

\par


\bibhang=1.7pc
\bibsep=2pt
\fontsize{9}{14pt plus.8pt minus .6pt}\selectfont
\renewcommand\bibname{\large \bf References}
\expandafter\ifx\csname
natexlab\endcsname\relax\def\natexlab#1{#1}\fi
\expandafter\ifx\csname url\endcsname\relax
  \def\url#1{\texttt{#1}}\fi
\expandafter\ifx\csname urlprefix\endcsname\relax\def\urlprefix{URL}\fi



\vskip .65cm
\noindent
Li-Xin Zhang

\vskip 2pt
\noindent School of Statistics and Data Sciences, Zhejiang Gongshang University, 310018

\vskip 2pt
\noindent
E-mail: stazlx@mail.zjgsu.edu.cn

\vskip 2pt
\noindent
ORCID ID: https://orcid.org/0000-0002-0121-0187

\newpage

 \centerline{\Large\bf Efficient Covariate-Adaptive
Randomization  }
\vspace{2pt} \centerline{\Large\bf  with Smallest Selection Bias in Clinical Trials }
\vspace{2pt}
 \centerline{\large\bf Li-Xin Zhang}
\vspace{.25cm}
 \author{Lixin Zhang$^{*}$}
\vspace{.4cm}
 \centerline{\it Zhejiang Gongshang University } 
\vspace{.55cm}
 \centerline{\Large\bf Supplementary Material}
\vspace{.55cm}
\fontsize{9}{11.5pt plus.8pt minus .6pt}\selectfont
\noindent
\par

\setcounter{section}{0}
\setcounter{equation}{0}
\def\theequation{S\arabic{section}.\arabic{equation}}
\def\thesection{S\arabic{section}}

\fontsize{12}{14pt plus.8pt minus .6pt}\selectfont

The supplementary material provides the theoretical proof of all the theorems, extended discussion on the selection bias, and the simulation study of the new CAR with a specific allocation function $g$. The results on the power are shown in Section \ref{sect:proof:power}. The properties of the CAR procedures are derived in Section  \ref{sect:proof:CAR}. The selection bias of the Pocock and Simon (1975) procedure  and its generalization  is  derived in Section \ref{sect:proof:SB}, while further discussion on the selection bias and randomness is given in Section \ref{sect:proof:ran}. The simulation study is presented  in Section  \ref{sect:simu}.  

\section{Deriving the asymptotic power}\label{sect:proof:power}
\setcounter{equation}{0}

In this section, we  prove the results
 on the power of   hypothesis testing to compare treatment effects which are given in  Section 2.

 \begin{center}\textbf{}
\end{center}
\vspace{-0.3in}
{\bf Proof of Theorem 1.} Write $W_n=(D_n/\sqrt{n})^2+Q_n^2$. Then
$$n-\ell_n^2=W_n+O(1)\frac{W_n^4}{n}=W_n\big(1+o(1)\big). $$
Let  $H(x)=\ln F(t_{1-\alpha}(\nu);\nu,\sqrt{x})$. Then
 $$H^{\prime}(x)=\frac{1}{2}\left.\frac{\frac{d}{d\delta}F(t_{1-\alpha}(\nu);\nu,\delta)}{\delta F(t_{1-\alpha}(\nu);\nu,\delta)}\right|_{\delta=\sqrt{x}}. $$
 We first show that
 \begin{equation}\label{eqproofPower.2}\frac{\frac{d}{d\delta}F(t;\nu,\delta)}{\delta F(t;\nu,\delta)}\to -1 \; \text{ as } \delta \to \infty \text{ and }\nu \to \infty
 \end{equation}
 uniformly in $t$ on a bounded interval. Notice that
 $$ F(t;\nu,\delta)=\ep\left[\Phi\big(t\sqrt{\chi^2(\nu)/\nu}-\delta\big)\right], $$
  $$ \frac{d}{d\delta} F(t;\nu,\delta)=-\ep\left[\varphi\big(t\sqrt{\chi^2(\nu)/\nu}-\delta\big)\right], $$
  where $\chi^2(\nu)$ is a random variable which has a $\chi^2$ distribution with degree of freedom $\nu $.
  Choose a sequence  $\epsilon_{\nu}$ such that $0<\epsilon_{\nu} \to 0$, $\epsilon_{\nu}\delta^2 \to \infty$ and $\epsilon_\nu \nu\to \infty$. Let $E=\{\chi^2(\nu)/\nu\le \epsilon_{\nu}\delta^2 \}$.
  Note
  $$ \varphi(\delta)\frac{\delta}{1+\delta^2}\le 1-\Phi(\delta)=\Phi(-\delta)\le \varphi(\delta)\frac{1}{\delta},\;\; \varphi(-\delta)=\varphi(\delta). $$
  On the event $E$,
  $$\frac{\varphi\big(t\sqrt{\chi^2(\nu)/\nu}-\delta\big)}{\delta \Phi\big(t\sqrt{\chi^2(\nu)/\nu}-\delta\big)}\to 1 $$
  uniformly in $t$ on a bounded interval and $\omega$. Hence
  $$ \ep\left[\varphi\big(t\sqrt{\chi^2(\nu)/\nu}-\delta\big)I_E\right]=\big(1+o(1)\big)\delta \ep\left[\Phi\big(t\sqrt{\chi^2(\nu)/\nu}-\delta\big)I_E\right]$$
  uniformly in $t$ on a bounded interval. On the other hand,
 \begin{align*}\pr(E^c)\le & e^{-\nu \epsilon_{\nu}\delta^2/4}\ep\left[e^{\chi^2(\nu)/4}\right]\le e^{-\nu (\epsilon_{\nu}\delta^2/4-1)}\le e^{-2\delta^2} \\
  \le & e^{-\delta^2}\delta \Phi(-\delta)\le e^{-\delta^2} \delta F(t;\nu,\delta)
  \end{align*}
  when $\delta$ and $\nu$ are large enough.
   Hence
  $$ \ep\left[\varphi\big(t\sqrt{\chi^2(\nu)/\nu}-\delta\big)\right]=\big(1+o(1)\big)\delta \ep\left[\Phi\big(t\sqrt{\chi^2(\nu)/\nu}-\delta\big)\right]$$
  uniformly in $t$ on a bounded interval. \eqref{eqproofPower.2} is proved, which implies
\begin{align*}H^{\prime}(x)=\frac{1}{2}\left.\frac{\frac{d}{d\delta}F(t_{1-\alpha}(\nu);\nu,\delta)}{\delta F(t_{1-\alpha}(\nu);\nu,\delta)}\right|_{\delta=\sqrt{x}}=-\frac{1}{2}\big(1+o(1))
\end{align*}
 as $\nu$ and $x$ tend to infinity. So, for $0\le y_n\le x_n\to \infty$ with $y_n/x_n\to 0$, we have
\begin{align*}
H(x_n^2)-H(x_n^2-y_n^2)=y_n^2H^{\prime}(z_n)=-\frac{y_n^2}{2}\left(1+o(1)\right),
\end{align*}
where $z_n\in (x_n^2-y_n^2,x_n^2)$.  Hence
\begin{equation}\label{eqPhiApp}\frac{F\left(t_{1-\alpha}(\nu); \nu, x_n\right)}{F\left(t_{1-\alpha}(\nu); \nu, \sqrt{x_n^2-y_n^2}\right)}=\exp\left\{-\frac{y_n^2}{2}\left(1+o(1)\right)\right\}.
\end{equation}
Now, let $x_n^2=\frac{\mu^2}{4\sigma_{\epsilon}^2}n$, $y_n^2=\frac{\mu^2}{4\sigma_{\epsilon}^2}\left(n- \ell_n^2\right)$. Then
$
y_n^2= \frac{\mu^2}{2\sigma_{\epsilon}^2}W_n^2(1+o(1)).
$
It follows that
\begin{equation}\label{eqproofPower.1} LossP_{T,n}(\mu|\bm X)=\exp\left\{ -\frac{u^2}{8\sigma_{\epsilon}^2}W_n\big(1+o(1)\big)\right\}
\end{equation}

Finally, note
$ \overline{\bm X}_n\to 0\;\; a.s. $
and
\begin{align*}
 \frac{S_{n,xx}}{n}=& \frac{1}{n} \left\{\sum_{i=1}^n (\bm X_i-\overline{\bm X}_n)^T(\bm X_i-\overline{\bm X}_n)\right.
 \\
 & \left.-N_{n,1}\big(\overline{\bm X}_{n,1}^T\overline{\bm X}_{n,1}-\overline{\bm X}_n^T\overline{\bm X}_n\big)-N_{n,2}\big(\overline{\bm X}_{n,2}^T\overline{\bm X}_{n,2}-\overline{\bm X}_n^T\overline{\bm X}_n\big)\right\}\\
 \to & \Var(\bm X)=diag(\sigma_{x,1}^2,\cdots,\sigma_{x,I}^2)\;\; a.s.
\end{align*}
We have
\begin{align*}
 W_n = & (D_n/\sqrt{n})^2+  \big(\bm D_n^x/\sqrt{n}\big)^T\big[\Var(\bm X)\big]^{-1}\big(\bm D_n^x/\sqrt{n}\big)\\
 & \qquad \quad +o(1)\left((D_n/\sqrt{n})^2+\big\|\bm D_n^x/\sqrt{n}\big\|^2\right) \\
 =& V_n \big(1+o(1)\big).
 \end{align*}
This completes the proof.  $\Box$

\bigskip
{\bf Proof of Theorem 2.} Notice that $\ep[X_k|d_k(X_k)]$ is a function of $\tilde{X_k}=d_k(X_k)$. We write it as $f_k(\tilde{X}_k)$. Then $X_{i,k}=\delta_{i,k}+f_k(\tilde{X}_{i,k})$. Under the condition that $\frac{D_n(k;t_k)}{\sqrt{n}}\to 0$, $t_k=1,\ldots, m_k$, $k=1,\ldots, p$, we have
$$ \frac{\sum_{i=1}^n (2T_i-1)f_k(\tilde{X}_{i,k})}{\sqrt{n}}=\frac{\sum_{t_k=1}^{m_k} D_n(k;t_k)f_k(\tilde{x}_k^{t_k})}{\sqrt{n}}\to 0. $$
Because the probability of $T_i=1$ is a function of $\tilde{\bm X}_1,\ldots,\tilde{\bm X}_i$,  the sequence $(2T_i-1)\big(\delta_{i,1},\ldots, \delta_{i,I}\big)$, $i=1,2,\ldots,$ is a sequence of  martingale vector differences with
$$ \Var\left((2T_i-1)\big(\delta_{i,1},\ldots, \delta_{i,I}\big)\big|\mathscr{F}_{i-1}\right)=diag\left(\sigma_{\delta,1}^2,\ldots,\sigma_{\delta,I}^2\right), $$
where $\mathscr{F}_i$ is the history $\sigma$-field generated by $T_1,\ldots, T_{i-1}$,  $\tilde{\bm X}_1,\ldots,\tilde{\bm X}_{i-1}$. By the central limit theorem of martingales, we have
$$\frac{\sum_{i=1}^n(2T_i-1)\big(\delta_{i,1},\ldots, \delta_{i,I}\big)^T}{\sqrt{n}}\overset{d}\to N_I\left(\bm 0,diag\left(\sigma_{\delta,1}^2,\ldots,\sigma_{\delta,I}^2\right)\right). $$
Hence
$$ V_n \overset{d}\to \sum_{k=1}^I \frac{\sigma_{\delta,k}^2}{\sigma_{x,k}^2}\chi_k^2(1). $$
The result follows. $\Box$

\bigskip
{\bf Proof of Remark 1.} We first consider the two-sided $t$-test. Notice that the distribution function of  $T^2$  is a F-distribution function $F(t;m,\nu,\delta)$ with degree of freedom $m=1$ and $\nu=n-p-2$, and the non-central parameter $\delta=\frac{\mu^2}{4\sigma_{\epsilon}^2}\ell_n^2$, it is sufficient to show that
\begin{equation}\label{eqFdistribution} \frac{\frac{d}{d\delta}F(t;m,\nu,\delta)}{F(t;m,\nu,\delta)}=-\frac{1}{2}\big(1+o(1)\big)
\end{equation}
uniformly in $t\in [0,t_0]$ as $\nu\to \infty$ and $\delta\to \infty$.  The non-central F-distribution has the following expansion
\begin{align*}
F(t;m,\nu,\delta)=&\sum_{j=0}^{\infty} e^{-\delta/2}\frac{(\delta/2)^j}{j!}F\left(\frac{m}{m+2j}t; m+2j,\nu,0\right)\\
=&\sum_{j=0}^{\infty} e^{-\delta/2}\frac{(\delta/2)^j}{j!}\int_0^tf(x;m,\nu,j)dx,
\end{align*}
where
$$f(x;m,\nu,j)=\frac{\Gamma\big((m+\nu)/2+j\big)}{\Gamma(m/2+j)\Gamma(\nu)}\frac{\nu^{\nu/2}m^{m/2+j}x^{m/2-1+j}}{(\nu+m x)^{(\nu+m)/2+j } }. $$
So,
\begin{align*}
&\frac{d F(t;m,\nu,\delta)}{d\delta} + \frac{1}{2} F(t;m,\nu,\delta)\\
&\; =\frac{1}{2}\sum_{j=0}^{\infty} e^{-\delta/2}\frac{(\delta/2)^j}{j!} \int_0^t f(x;m,\nu,j+1)dx\\
& \; =   \frac{1}{2}\sum_{j=0}^{\infty} e^{-\delta/2}\frac{(\delta/2)^j}{j!}F\left(\frac{m}{m+2(j+1)}t; m+2(j+1),\nu,0\right).
\end{align*}
It is obvious that
\begin{align*}
& e^{-\delta/2}\frac{(\delta/2)^j}{j!}F\left(\frac{m}{m+2(j+1)}t; m+2(j+1),\nu,0\right)\\
=  & \frac{2(j+1)}{\delta}e^{-\delta/2}\frac{(\delta/2)^{j+1}}{(j+1)!}F\left(\frac{m}{m+2(j+1)}t; m+2(j+1),\nu,0\right)\\
\le  &\frac{2(j+1)}{\delta}  F(t;m,\nu,\delta)
\end{align*}
and
$$ \int_0^tf(x;m,\nu,j+1)dx\le \frac{\nu+m+2j}{m+2j}\frac{mt}{\nu+mt} \int_0^tf(x;m,\nu,j)dx. $$
It follows that
\begin{align*}
0\le & \frac{d F(t;m,\nu,\delta)}{d\delta} + \frac{1}{2} F(t;m,\nu,\delta) \\
\le & \frac{1}{2}\sum_{j=0}^{Km}\frac{2(j+1)}{\delta}  F(t;m,\nu,\delta)\\
& +\frac{1}{2}\sum_{j=Km}^{\infty} e^{-\delta/2}\frac{(\delta/2)^j}{j!}  \left(\frac{mt}{\nu}+\frac{mt}{m+2j}\right)\int_0^tf(x;m,\nu,j)dx\\
\le & \left(\frac{2(Km+1)^2}{\delta}+\frac{mt}{\nu}+\frac{t}{2K+1}\right)F(t;m,\nu,\delta).
\end{align*}
This proves \eqref{eqFdistribution}.

For the one-sided $z$-test for 
$$
H_0: \mu=0  \mbox{ versus } H_A: \mu > 0.
$$
where $\mu=\mu_1-\mu_2$ is the treatment effect difference, the test statistic is
\begin{equation}\label{z:stat}
U=\frac{\bm L \boldsymbol{\hat\beta}}{(\sigma_{\epsilon}^2\bm L(\boldsymbol{X}^T\boldsymbol{X})^{-1}\bm L^T)^{1/2}}.
\end{equation}
The power function is
$$\beta_{U,n}(\mu|\bm X)=1-\Phi(z_{1-\alpha}-\sqrt{\delta}) \; \text{ with } \delta =\frac{u^2}{4\sigma_{\epsilon}^2}\ell_n^2 . $$
It is obvious that
$$ \frac{\frac{d}{d\delta}\Phi(z_{1-\alpha}-\sqrt{\delta})}{ \Phi(z_{1-\alpha}-\sqrt{\delta})}\to -\frac{1}{2}\text{ as } \delta\to \infty. $$
We obtain a  result similar to \eqref{eqFdistribution}.

For the two-sided $z$-test for
$$
H_0: \mu=0  \mbox{ versus } H_A: \mu \ne  0,
$$
the test statistic is $U$ and the power function is
$$\beta_{|U|,n}(\mu|\bm X)=1-\Phi(z_{1-\alpha/2}-\sqrt{\delta})+\Phi(-z_{1-\alpha/2}-\sqrt{\delta}) \; \text{ with } \delta =\frac{u^2}{4\sigma_{\epsilon}^2}\ell_n^2 . $$
It is easily seen that
\begin{align*} & \frac{\frac{d}{d\delta}\left[\Phi(z_{1-\alpha/2}-\sqrt{\delta})-\Phi(-z_{1-\alpha/2}-\sqrt{\delta})\right]}{\Phi(z_{1-\alpha/2}-\sqrt{\delta})-\Phi(-z_{1-\alpha/2}-\sqrt{\delta})}\\
=&-\frac{1}{2\sqrt{\delta}}\frac{ \varphi(z_{1-\alpha/2}-\sqrt{\delta})-\varphi(-z_{1-\alpha/2}-\sqrt{\delta})}{\Phi(z_{1-\alpha/2}-\sqrt{\delta})-\Phi(-z_{1-\alpha/2}-\sqrt{\delta})}\\
\sim & -\frac{1}{2\sqrt{\delta}}\frac{ \varphi(z_{1-\alpha/2}-\sqrt{\delta})}{\Phi(z_{1-\alpha/2}-\sqrt{\delta})}\to -\frac{1}{2}\;\; \text{ as }\delta \to \infty.
\end{align*}
We obtain a results   similar to \eqref{eqFdistribution}. This completes the proof. $\Box$

\bigskip

\section{Deriving the properties of the CAR procedures}\label{sect:proof:CAR}
\setcounter{equation}{0}

In this section, we  prove the results
 on covariate adaptive randomization procedures which are given in Section 4.

Recall  that $\bm t= (t_1,\ldots,t_I)$, $t_i=1,\ldots,m_i, i=1,\ldots, I$, has $M=\prod_{k=1}^{I}m_k$ values.
Our purpose is to  study the properties of $\bm D_n=[D_n(\bm t)]$. Besides $\bm D_n$, we will also consider the weighted average of the imbalances $\Lambda_{n-1}(\bm t)$ as in (2.1). Let
 \begin{align*}
\Lambda_{n}(\bm t)=& w_{o}D_{n}+\sum_{i=1}^{I}w_{m,i}D_{n}(i;t_i)+w_{s}D_n(\bm t),
\\
 \bm \Lambda_{n}=& \left[\Lambda_{n}(\bm t)\right]_{1\leq t_1\leq m_1,\ldots, 1\leq t_I\leq m_I}.
\end{align*}
Also let $\mathscr{F}_{n-1}$ be the history $\sigma$-field generated by  the covariates $\tilde{\bm X}_1,\ldots, \tilde{\bm X}_{n-1}$ and   results of allocation $T_1,\ldots, T_{n-1}$. Then
the allocation probability  of the $n$th patient is
$$ \pr\left(T_n=1\big|\mathscr{F}_{n-1}, \tilde{\bm X}_n=(\tilde{x}_1^{t_1},\ldots, \tilde{x}_I^{t_I})\right)=g_n\left(4\Lambda_{n-1}(\bm t)\right). $$
 is a function of $\bm \Lambda_{n-1}$.      It is obvious that $\bm \Lambda_n=\bm L(\bm D_n):\bm D_n \to \bm \Lambda_n$ is a linear transform of $\bm D_n$.
The following proposition gives the relation between $\bm D_n$ and $\bm\Lambda_n$ and tells us that both $(\bm D_n)_{n\ge 1}$ and $(\bm\Lambda_n)_{n\ge 1}$ are Markov chains.
\begin{proposition}\label{proposition1}  (i) If $w_s>0$, then $\bm \Lambda_n=\bm L(\bm D_n)$ is a one to one linear map; If $w_s+w_{m,i}>0$, then each $D_n(i;t_i)=D_{i;t_i}(\bm\Lambda_n)$ is a linear transform of $\bm \Lambda_n$; For any case, $D_n=D(\bm\Lambda_n)$ is a linear transform of $\bm \Lambda_n$.

(ii) $(\bm{D}_{n})_{n\geq 1}$ is a  non-homogeneous Markov chain on the space $\mathbb{Z}^{m}$ with period 2;

(iii)  $(\bm \Lambda_n)_{n\geq 1}$ is a  non-homogeneous Markov chain on the space $\bm L(\mathbb{Z}^{m})$ with period 2.
\end{proposition}

\smallskip

\noindent{\bf Proof.} For (i), taking the summation of $\Lambda_n(\bm t)$ over all $\bm t$ yields
$$ \sum_{\bm t} \Lambda_n(\bm t)=\big(w_s + \sum_{i=1}^I w_{m,i}\prod_{j\ne i} m_j  + w_o M\big) D_n. $$
So $D_n$ is a linear transform of $\bm \Lambda_n$.
Taking the summation of $\Lambda_n(\bm t)$ over all $t_1,\ldots, t_I$ except $t_i$  yields
\begin{align*}
 &\sum_{t_1,\ldots,t_{i-1}, t_{i+1},\ldots, t_I} \Lambda_n(t_1,\ldots,t_I)\\
 =&\big(w_s + w_{m,i}\prod_{j\ne i} m_j)D_n(i;t_i)+\big( \sum_{l\ne i} w_{m,l}\prod_{j\ne i, l} m_j  + w_o \prod_{j\ne i} m_j\big) D_n.
 \end{align*}
Hence, when $w_s+w_{m,i}>0$, each $D_n(i;t_i)$ is a linear transform of $\bm \Lambda_n $ and $D_n$, and so it is  a linear transform of $\bm \Lambda_n$. Finally, when $w_s>0$, it is obvious that each $D_n(\bm t)$ is a linear transform of $\Lambda_n(\bm t)$, $D_n(1;t_1),\ldots, D_n(I;t_I)$ and $D_n$, and so it is a linear transform of $\bm \Lambda_n$. Hence,  when $w_s>0$, $\bm \Lambda_n=\bm L(\bm D_n)$ is a one to one linear map.

\smallskip
For (ii) and (iii), notice
$$ D_n(\bm t)=D_{n-1}(\bm t)+2\big(T_n-\frac{1}{2}\big)\mathbb{I}\{\tilde{\bm X}_n=(\tilde{x}_1^{t_1},\ldots, \tilde{x}_I^{t_I})\}. $$
 Then
\begin{align*}
 \pr(\Delta D_n(\bm t)=1|\mathscr{F}_{n-1})= & g_n(4\Lambda_{n-1}(\bm t))p(\bm t), \\
 \pr(\Delta D_n(\bm t)=-1|\mathscr{F}_{n-1})= & \big[1-g_n\big(4\Lambda_{n-1}(\bm t)\big)\big]p(\bm t),\\
\pr(\Delta D_n(\bm t)=0|\mathscr{F}_{n-1})= &1-p(\bm t),
\end{align*}
where $p(\bm t)=p(t_1,\ldots,t_I)=\pr(\tilde{\bm X}_n=(\tilde{x}_1^{t_1},\ldots, \tilde{x}_I^{t_I}))$, $\bm p=\big(p(\bm t),\bm t\in \mathbb{Z}^{m_1\times \ldots\times m_I}\big)$.  Let $\Delta\mathscr{D}$ be the state space of $\Delta \bm D_n=\bm D_n-\bm D_{n-1}$, i.e., each $\bm d\in \Delta\mathscr{D}$ has only one non-zero element which is $1$ or $-1$.  For two vectors $\bm x$ and $\bm y$ on $\mathbb{Z}^{m_1\times \ldots\times m_I}$, we write $\bm x\cdot \bm y=\sum_{\bm t}x(\bm t)y(\bm t)$, $|\bm x|=(|x(\bm t)|)$. The conditional probability above can be written in the following form,
 \begin{align}\label{eqprobabDD}
  \pr(\Delta \bm D_n=\bm d|\mathscr{F}_{n-1})=&\left(g_n(4|\bm d|\cdot \bm \Lambda_{n-1})-\frac{1}{2}\right)\bm d\cdot \bm p+\frac{1}{2}|\bm d\cdot \bm p|\\
  =& \left(g_n\left(4|\bm d|\cdot \bm L(\bm D_{n-1})\right)-\frac{1}{2}\right)\bm d\cdot \bm p+\frac{1}{2}|\bm d\cdot \bm p|, \nonumber\\
  &\;\; \bm d\in \Delta \mathscr{D},\nonumber
  \end{align}
 which depends only on $\bm \Lambda_{n-1}=\bm L(\bm D_{n-1})$. So, conditional on $\bm{D}_{n-1}$, $\bm D_n$ is  conditionally independent of $(\bm{D}_1,\ldots,\bm{D}_{n-2})$. It follows that $(\bm D_n)_{n\ge 1}$ is a Markov chain on $\mathbb{Z}^{m}$.

For (iii), notice that for any point $\bm e$ in the state space $\{ \bm L(\bm d):\bm d\in \Delta\mathscr{D}\}$ of $\Delta\bm \Lambda_n$,
\begin{align}\label{tran:probab}
 & \pr(\Delta \bm \Lambda_n=\bm e|\mathscr{F}_{n-1})=\nonumber\\
  &\sum_{\bm d:  \bm L(\bm d) =\bm e,\bm d\in \Delta\mathscr{D}} \left\{\left(g_n(4|\bm d|\cdot \bm \Lambda_{n-1})-\frac{1}{2}\right)\bm d\cdot \bm p+\frac{1}{2}|\bm d\cdot \bm p|\right\},
  \end{align}
which depends only on $\bm \Lambda_{n-1}$. So, given $\bm \Lambda_{n-1}$, $\bm \Lambda_n$ is  conditionally independent of $(\bm \Lambda_1,\ldots,\bm \Lambda_{n-2})$. It follows that
$\bm \Lambda_n$ is a Markov chain. 
  This completes the proof of Proposition \ref{proposition1}.
$\Box$

\bigskip

\noindent {\bf Proof of Theorem 4.}
Recall
 $$ M_n=\sum_{\bm t}w_s D_n^2(\bm t)+\sum_{i=1}^I\sum_{t_i=1}^{m_i} w_{m,i} D_n^2(i;t_i)+w_o D_n^2. $$
By Proposition \ref{proposition1} (i), $M_n=M^{\ast}(\bm\Lambda_n)$ is a   function of $\bm \Lambda_n$, and
$ M_n\le C\|\bm \Lambda_n\|^2. $
On the other hand,
 \begin{align*}
& |\Lambda_{n}(\bm t)|^2\le \big(w_{o}|D_{n}|+\sum_{i=1}^{I}w_{m,i}|D_{n}(i;t_i)|+w_{s}|D_n(\bm t)|\big)^2\\
& \le  \big(w_{o}|D_{n}|^2+\sum_{i=1}^{I}w_{m,i}|D_{n}(i;t_i)|^2+w_{s}|D_{n}(\bm t)|^2\big)(w_{o} +\sum_{i=1}^{I}w_{m,i} +w_{s} ),\\
&=  w_{o}|D_{n}|^2+\sum_{i=1}^{I}w_{m,i}|D_{n}(i;t_i)|^2+w_{s}|D_{n}(\bm t)|^2,
\end{align*}
which implies that $ \|\bm \Lambda_n\|^2\le M\cdot M_n$.
We will prove the theorem via two steps. First, we will show that under condition 
\begin{equation}
    g_n(x)\le 1/2 \le g_n(-x)  \text{ when } x\ge 0, \label{eqCondF2}
\end{equation}
we have
\begin{equation}\label{eqLrBound}
M_n=O(n) \;\;\text{in } L_r\;\; \forall r>0.
 \end{equation}
 In the second step, we will show that under conditions \eqref{eqCondF2} and 
 \begin{equation}
 \frac{\left|1/2-g_n(x_n)\right|}{|x_n|/n}\to +\infty \; \text{ whenever }\; 0<\frac{|x_n|}{n}\to 0, \label{eqCondF3}
 \end{equation}
 we have
 \begin{equation}\label{eqConvInProbab}
M_n=o(n) \;\;\text{ in probability},
 \end{equation}
 which, together with \eqref{eqLrBound}, implies that $M_n=o(n)$ in $L_r$ for any $r>0$.
 (i)-(iii) follow from  immediately from this. Finally,
\begin{align*}
\ep\left[p_n\vee(1-p_n)\right]=&\frac{1}{2}+ \sum_{\bm t}
p(\bm t)\ep\Big| \frac{1}{2}-g_n\big(4 \Lambda_{n-1}(\bm t) \big)\Big|
  \to \frac{1}{2}
\end{align*}
due to the condition that
\begin{equation}\label{eqCondF4}
  g_n(x_n)\to \frac{1}{2} \text{ whenever }\; \frac{x_n}{n}\to 0,
  \end{equation}
and the fact that $\frac{\|\bm \Lambda_n\|}{n}\le C\frac{M_n^{1/2}}{n}\to 0$ in probability. Hence
$$ SB_n=\frac{1}{n}\sum_{m=1}^n \ep\left[p_m\vee(1-p_m)\right]\to \frac{1}{2}. $$

\smallskip
Now, we begin the proofs of \eqref{eqLrBound} and \eqref{eqConvInProbab}.   Given $\tilde{\bm X}_n=(\tilde{x}_1^{t_1},\ldots, \tilde{x}_I^{t_I})$, if $T_n=1$, then
 \begin{align*}
 M_n-M_{n-1}=& w_s \left\{ \big(D_{n-1}(\bm t)+1\big)^2-D_{n-1}^2(\bm t)\right\}\\
 &+\sum_{i=1}^I w_{m,i} \left\{ \big(D_{n-1}(i;t_i)+1\big)^2-D_{n-1}^2(i;t_i)\right\}\\
 &+  w_o \left\{ (D_{n-1}+1)^2-D_{n-1}^2 \right\}\\
 =& 2\Lambda_{n-1}(\bm t)+1,
 \end{align*}
 while, if $T_n=0$, then
 $M_n-M_{n-1}= -2\Lambda_{n-1}(\bm t)+1. $
 So,
\begin{equation}\label{eqV} M_n-M_{n-1}=2\Lambda_{n-1}(\bm t)\left(2T_n-1\right)\mathbb{I}\{\tilde{\bm X}_n=(\tilde{x}_1^{t_1},\ldots, \tilde{x}_I^{t_I})\}+1.
\end{equation}
 Hence
\begin{align*}
  &\ep\left[M_n-M_{n-1}\big|\tilde{\bm X}_n=(x_1^{t_1},\ldots, x_I^{t_I}),\mathscr{F}_{n-1}\right]\\
 =&2\Lambda_{n-1}(\bm t)\left[ 2g_n(4\Lambda_{n-1}(\bm t))-1\right]+1\\
 =& -4\Big|\Lambda_{n-1}(\bm t)\Big|\cdot\left| \frac{1}{2}-g_n\left(4\Lambda_{n-1}(\bm t)\right)\right|+1,
 \end{align*}
due to the condition that $g_n(-x)\le 1/2\le g_n(x)$ when $x\ge 0$. It follows that
 \begin{align}\label{eqcondE}
  \ep\left[M_n\big|\mathscr{F}_{n-1}\right]-M_{n-1}
  = -4S_n(\bm\Lambda_{n-1})+1,
 \end{align}
 where
  \begin{equation}\label{Sfunction}
  S_n(\bm\Lambda_{n-1})=\sum_{\bm t}\Big|\Lambda_{n-1}(\bm t)\Big|\cdot\left| \frac{1}{2}-g_n\left(4\Lambda_{n-1}(\bm t)\right)\right|p(\bm t)
  \end{equation}
 is a nonnegative function of $\bm \Lambda_{n-1}$.

From \eqref{eqcondE} and by noting $M_0=0$, it follows that
 \begin{equation}\label{eq2edM} \ep\left[M_n\right]\le n.
 \end{equation}

 Further,     by \eqref{eqV} we have
\begin{align*}
 M_n=&M_{n-1}+1+2\Lambda_{n-1}(\bm t)(2T_n-1)\mathbb{I}\{\tilde{\bm X}_n=(\tilde{x}_1^{t_1},\ldots, \tilde{x}_I^{t_I})\}\\
 \widehat{=}& M_{n-1}+1+\xi .
 \end{align*}
It is obvious that
$$|\xi|=2|\Lambda_{n-1}(\bm t) |\le 2 \sqrt{ M_{n-1}},\;\; \ep[\xi|\mathscr{F}_{n-1}]=-4S_n(\bm\Lambda_{n-1}). $$
It follows that for positive integer $r$,
\begin{align*}
M_n^{r+1}- & M_{n-1}^{r+1}=(r+1)(M_{n-1}+1)^r\xi\\
&\quad +\left\{(M_{n-1}+1)^{r+1}-M_{n-1}^{r+1}+\sum_{k=2}^{r+1}\binom{r+1}{k}\xi^k(M_{n-1}+1)^{r+1-k}\right\}\\
&\qquad  \le (r+1)(M_{n-1}+1)^r\xi+C_r (M_{n-1}+1)^r,
\end{align*}
where $C_r$ is a constant which depends on $r$. It follows that
\begin{align} \label{eq:conditionalM}
&\ep[M_n^{r+1}|\mathscr{F}_{n-1}]-M_{n-1}^{r+1}\nonumber\\
\le & -4 (r+1)(M_{n-1}+1)^rS_n(\bm\Lambda_{n-1})+C_r (M_{n-1}+1)^r\\
\le &C_r (M_{n-1}+1)^r.\nonumber
\end{align}

By \eqref{eq2edM} and the induction we conclude \eqref{eqLrBound}. Obviously, \eqref{eqLrBound} implies
$$ \frac{\|\bm \Lambda_n\|}{n}\to 0 \text{ in probability}. $$

Next, we consider \eqref{eqConvInProbab}. By \eqref{eqcondE},
 \begin{align}\label{eqexpecationofV}
  \ep\left[M_n\right]-\ep\left[M_{n-1}\right]
  = -4\ep\left[S_n(\bm\Lambda_{n-1})\right]+1.
 \end{align}
 Let $m:=m_n\in\{0,1,\ldots, n\}$ be the last one for which $1-4\ep\left[S_{m+1}(\bm\Lambda_{m})\right]\ge 0$. Then
 $$ \ep\left[M_n\right] \le \ep\left[M_{m+1}\right]\le \ep\left[M_m\right]+1 $$
 and
 $$4\ep\left[S_{m+1}(\bm\Lambda_{m})\right]\le 1.  $$
 If $m=m_n$ is bounded, the proof is completed. Assume $m\to \infty$. Recall \eqref{Sfunction}.  Notice that the condition \eqref{eqCondF3} implies
 $$ \frac{S_{m+1}(\bm\Lambda_{m})}{\frac{1}{m}\sum_{\bm t} |\Lambda_m(\bm t)|^2p(\bm t)}\to +\infty \;\; as\; \frac{\|\bm \Lambda_m\|}{m}\to 0. $$
 On the other hand,
  $$M_m\ge \sum_{\bm t} |\Lambda_m(\bm t)|^2p(\bm t)\ge c_0  M_m. $$
  So for any $0<\epsilon<1$, there exists a $\delta>0$ such that
 $$ S_{m+1}(\bm\Lambda_{m})I\Big\{\frac{\|\bm \Lambda_m\|}{m}\le \delta\Big\}\ge \frac{1}{\epsilon}  \frac{1}{m} M_mI\Big\{\frac{\|\bm \Lambda_m\|}{m}\le \delta\Big\}. $$
 So
 $$\ep\Big[M_mI\Big\{\frac{\|\bm \Lambda_m\|}{m}\le \delta\Big\}\Big]\le \epsilon m \ep\Big[S_{m+1}(\bm\Lambda_{m})\Big]\le   \frac{\epsilon m}{4}. $$
 Also,
\begin{align*} \ep\Big[M_mI\Big\{\frac{\|\bm \Lambda_m\|}{m}> \delta\Big\}\Big]\le & \Big\{\ep\Big[M_m^2\Big]\Big\}^{1/2}\Big\{\pr\Big(\frac{\|\bm \Lambda_m\|}{m}> \delta\Big)\Big\}^{1/2}\\
\le &c m \pr^{1/2}\Big(\frac{\|\bm \Lambda_m\|}{m}> \delta\Big).
\end{align*}
It follows that
$$ \ep\left[M_n\right]\le n\Big\{\frac{\epsilon}{4}+c\pr^{1/2}\Big(\frac{\|\bm \Lambda_m\|}{m}> \delta\Big)\Big\}+1. $$
By the fact  $\pr\Big(\frac{\|\bm \Lambda_m\|}{m}> \delta\Big)\to 0$, we conclude that
$\frac{\ep\left[M_n\right]}{n}\to 0. $
This proves \eqref{eqConvInProbab}. The proof of Theorem 4 is completed. $\Box$

 \bigskip
 \noindent {\bf Proof of Theorem 5.}
 Notice that, by \eqref{eq:conditionalM},
 \begin{align*} \ep[M_n^{r+1}]-\ep[M_{n-1}^{r+1}]\le & -4 (r+1)\ep\left[(M_{n-1}+1)^rS_n(\bm\Lambda_{n-1})\right]+C_r \ep\left[(M_{n-1}+1)^r\right].
\end{align*}
Let $m:=m_n\in\{0,1,\ldots, n\}$ be the last one for which
\begin{equation}\label{eqinduction1}-4 (r+1)\ep\left[(M_m+1)^rS_{m+1}(\bm\Lambda_m)\right]+C_r \ep\left[(M_m+1)^r\right]\ge 0. \end{equation}
Then
\begin{equation}\label{eqinduction2} \ep[M_n^{r+1}]\le \ep[M_{m+1}^{r+1}]\le \ep[M_m^{r+1}]+C_r \ep\left[(M_m+1)^r\right].\end{equation}
Recall \eqref{Sfunction}. Note for $\delta>0$ small enough, on the event $E=E(\bm t)=\left\{\frac{|\Lambda_{m}(\bm t)|}{m^{\gamma}}\le \delta\right\}$,
$$  S_{m+1}(\bm\Lambda_m)\ge \frac{-g^{\prime}(0)}{2}\frac{1}{m^{\gamma}}|\Lambda_m(\bm t)|^2 p(\bm t)\ge c_0 \frac{1}{m^{\gamma}}|\Lambda_m(\bm t)|^2, $$
and on the event $E^c$,
$$  S_{m+1}(\bm\Lambda_m)\ge \big(\frac{1}{2}-g(c\delta) \big)\vee \big(\frac{1}{2}-g(-c\delta) \big)  |\Lambda_m(\bm t)|  p(\bm t)\ge c_0 |\Lambda_m(\bm t)|. $$
 $ C^{-1} M_n\le \|\bm \Lambda_n\|^2\le C\cdot M_n$. Then by \eqref{eqinduction1},
$$ \ep[|\Lambda_m(\bm t)|^{2r+2}E]\le  C    m^{\gamma} \ep\left[(M_m+1)^rS_{m+1}(\bm\Lambda_m)\right] \le C_r   m^{\gamma} \ep\left[(M_m+1)^r\right], $$
and
$$\ep[|\Lambda_m(\bm t)|^{2r+1}E^c]\le  C     \ep\left[(M_m+1)^rS_{m+1}(\bm\Lambda_m)\right] \le C_r   \ep\left[(M_m+1)^r\right]. $$
By the H\"older inequality,
\begin{align*}
 &\ep[|\Lambda_m(\bm t)|^{2r+2}E^c]= \ep\left[|\Lambda_m(\bm t)|^{\frac{2r+1}{p}}|\Lambda_m(\bm t)|^{\frac{2r+1}{q}+1}E^C\right]\\
  \le &\left(\ep\left[|\Lambda_m(\bm t)|^{2r+1}E^c\right]\right)^{1/p}\left(\ep\left[|\Lambda_m(\bm t)|^{2r+1+q}E^c\right]\right)^{1/q}\\
 \le & C_r\left(\ep\left[(M_m+1)^r\right]\right)^{1/p}\left(\ep\left[M_m^{r+\frac{1}{2}+\frac{1}{2}q}\right]\right)^{1/q}\\
  \le & C_r\left(\ep\left[(M_m+1)^r\right]\right)^{1/p}\left(\ep\left[M_m^{2r+1+ q}\right]\right)^{1/2q},
 \end{align*}
 where $p,q>1$, $1/p+1/q=1$. It follows that
 \begin{align}\label{eqinductionsad}
   \ep[M_n^{r+1}]\le  & \ep[M_m^{r+1}]+C_r \ep\left[(M_m+1)^r\right]\le   c_r\sum_{\bm t}\ep[|\Lambda_m(\bm t)|^{2r+2}] +C_r \ep\left[(M_m+1)^r\right]\nonumber\\
 \le & C_r   m^{\gamma}\ep\left[(M_m+1)^r\right]+C_r\left(\ep\left[(M_m+1)^r\right]\right)^{1/p}\left(\ep\left[M_m^{2r+1+ q}\right]\right)^{1/2q}.
 \end{align}
 We will use the induction method to obtain the moment estimates of $\ep[M_n^r]$ by inducting both in $n$ and $r$.  We assume
 \begin{equation}\label{eqinduction3} \ep[M_n^s]\le C_s n^{\beta s} \text{ for all integers } n,s\ge 0,
 \end{equation}
 and
\begin{equation}\label{eqinduction4} \ep[M_n^r]\le C_r n^{\alpha r} \text{ for all } n,
\end{equation}
where $0<\beta/2<\alpha\le \beta\le 2$.

Let  $q$ be an integer large enough such that
$$ \frac{\alpha r}{p}+\beta\frac{2r+1}{2q}+\frac{1}{2}\beta=\alpha r +\frac{r(\beta-\alpha)+\beta/2}{q}+\frac{1}{2}\beta\le (r+1)\alpha. $$
Then
  \begin{align*}
 \ep[M_m^{r+1}]\le & C_r m^{\gamma} m^{\alpha r}+ C_r \left(m^{\alpha r}\right)^{1/p} \left(m^{\beta(2r+1+q)}\right)^{1/2q}\\
 \le & C_r\left(m^{\gamma} m^{\alpha r}+m^{\alpha(r+1)}\right).
 \end{align*}
 By \eqref{eqinduction2},
 \begin{equation}\label{eqinduction5} \ep[M_n^{r+1}]\le C_r\big(n^{\gamma} n^{\alpha r}+n^{\alpha(r+1)}\big)\le C_{r+1}n^{\alpha(r+1)}\;\; \text{ for all } n,
 \end{equation}
 whenever $\alpha\ge \gamma$.

 Obviously, \eqref{eqinduction3} holds for $\beta=2$ since $M_n\le n^2$,  and  \eqref{eqinduction4} is   true for $r=0$. So, by \eqref{eqinduction5} and the induction, \eqref{eqinduction4} holds for all integer $r$ whenever $\alpha>1$, which implies   that \eqref{eqinduction3} holds for any $\beta$ with   $\beta\ge \gamma$ and $\beta>1$. Repeating  the above procedure $i$ times concludes that \eqref{eqinduction3} holds for any $\beta$ with   $\beta\ge \gamma$ and $\beta>1/2^{i-1}$. Stopping the procedure when $1/2^{i-1}<\gamma$, we conclude that
 $$ \ep M_n^r\le C_r n^{\gamma r} \;   \text{ for all integers } n \text{ and } r. $$
 Thus,
 \begin{equation}\label{eqth2.1}
M_n=O(n^{\gamma}) \;\; \text{ in } L_r   \text{ for all } r>0,
\end{equation}
is proved. By  \eqref{eqth2.1}, for any $\epsilon,\delta>0$,
 $$ \sum_n \pr\left( M_n\ge \delta n^{\gamma +\epsilon} \right)\le \sum_n  \frac{\ep  M_n^Q}{\delta^Q n^{Q\gamma}n^{\epsilon Q}}
 \le C\le \sum_n  n^{-\epsilon Q}<\infty
 $$
 when $Q>1/\epsilon$. By the Borel-Cantelli lemma, 
\begin{equation}\label{eqth2.2}
M_n=o(n^{\gamma+\epsilon}) \; a.s.    \text{ for any } \epsilon>0
\end{equation}
 is proved. And so, (i)-(iii) follow from Proposition \ref{proposition1}.

 Finally, we consider the selection bias. Note
 \begin{align}\label{TaylaExp1}
 \left|\frac{1}{2}-g_{m+1}\left(4 \Lambda_m(\bm t) \right)\right|
 =-4g^{\prime}(0)\frac{\big|\Lambda_m(\bm t)\big|}{m^{\gamma}}+4\theta_m(\bm t) \frac{\big|\Lambda_m(\bm t)\big|}{m^{\gamma}},
 \end{align}
 where
 $$\theta_m(\bm t)=\frac{\left|\frac{1}{2}-g\left(4\frac{\Lambda_m(\bm t)}{m^{\gamma}}\right)\right|}{\frac{4|\Lambda_m(\bm t)|}{m^{\gamma}}}+4g^{\prime}(0)$$ is bounded and converges to 0  in probability due to the fact that $0\le g(x)\le 1$, $g^{\prime}(0)$ is finite  and $\frac{|\Lambda_m(\bm t)|}{m^{\gamma}}\to 0 $ in probability. So,
 $\ep\theta_m(\bm t)^2 =o(1) $ as $m\to \infty$.  By the Cauchy-Schwartz  inequality,
 $$ \ep\left|\theta_m(\bm t) \frac{\big|\Lambda_m(\bm t)\big|}{m^{\gamma}}\right|=o(1)\cdot   \left(\ep\left|\frac{\big|\Lambda_m(\bm t)\big|}{m^{\gamma}}\right|^2\right)^{1/2}. $$

It follows that
\begin{align}\label{TaylaExp2}
&\ep\left[p_m\vee(1-p_m)\right]= \frac{1}{2}+ \sum_{\bm t}
p(\bm t)\ep\left| \frac{1}{2}-g_m\left(4\Lambda_{m-1}(\bm t)\right)\right|  \nonumber \\
= & \frac{1}{2}+ \sum_{\bm t}
p(\bm t)\ep\left[- 4 g^{\prime}(0)\frac{|\Lambda_{m-1}(\bm t)|}{(m-1)^{\gamma}}\right]
 +\sum_{\bm t}
p(\bm t) o(1) \left(\ep\left|\frac{\big|\Lambda_{m-1}(\bm t)\big|}{{(m-1)}^{\gamma}}\right|^2\right)^{1/2}\nonumber\\
=&\frac{1}{2}- 4 g^{\prime}(0) \sum_{\bm t}
p(\bm t)\frac{\ep\left[|\Lambda_{m-1}(\bm t)|\right]}{(m-1)^{\gamma}}+o\left(\frac{1}{(m-1)^{\gamma/2}}\right),
\end{align}
where,   the last equality is due to the fact that $\|\bm \Lambda_m\|\le C M_m^{1/2}=O(m^{\gamma/2})$ in any $L_r$. Hence
\begin{align*}
 SB_n=&\frac{1}{n}\sum_{m=1}^n \ep\left[p_m\vee(1-p_m)\right]\\
 =& \frac{1}{2}-4 g^{\prime}(0)\frac{1}{n}\sum_{m=1}^{n} \sum_{\bm t}
p(\bm t)\frac{\ep\left[|\Lambda_{m}(\bm t)|\right]}{m^{\gamma}}+o\left(\frac{1}{n^{\gamma/2}}\right) \\
 =: &\frac{1}{2}+a_nn^{-\gamma/2}+o\left(n^{-\gamma/2}\right).
 \end{align*}

 Next, we show $a_n$ is bounded away from zero and infinity. From \eqref{eqexpecationofV}, we conclude that
 \begin{align*}
 0=\lim_{n\to \infty} \frac{\ep[M_{n+1}]}{n}=-4\lim_{n\to \infty}\frac{1}{n} \sum_{m=1}^n \ep[S_{m+1}(\bm \Lambda_m)]+1.
 \end{align*}
 By noting \eqref{TaylaExp1}, similar to \eqref{TaylaExp2} we have
 \begin{align*}\ep[S_{m+1}(\bm \Lambda_m)]=&-g^{\prime}(0) \sum_{\bm t} p(\bm t) \frac{\ep|\Lambda_m(\bm t)|^2}{m^{\gamma}}+o(1)\sum_{\bm t} p(\bm t) \ep\left(\frac{\ep|\Lambda_m(\bm t)|^2}{m^{ \gamma}}\right)^2\\
 =&-g^{\prime}(0) \sum_{\bm t} p(\bm t) \frac{\ep|\Lambda_m(\bm t)|^2}{m^{\gamma}}+o(1).
 \end{align*}
It follows that
  $$ \lim_{n\to \infty} \frac{1}{n}\sum_{m=1}^n \sum_{\bm t} p(\bm t) \frac{\ep|\Lambda_m(\bm t)|^2}{m^{\gamma}}=\frac{1}{-16g^{\prime}(0)}. $$
  Regarding the summation on the  right-hand side above as $\ep|Z_n|^2$ of a random variable $Z_n=\frac{\Lambda_{\tau}(\bm\kappa)}{\tau^{\gamma/2}}$, where $\pr(\tau=m)=\frac{1}{n}$, $\pr(\bm\kappa=\bm t)=p(\bm t)$, then
  $$\frac{1}{n}\sum_{m=1}^{n} \sum_{\bm t}
p(\bm t)\frac{\ep |\Lambda_{m}(\bm t)|}{m^{\gamma/2}}=\ep|Z_n|\le (\ep|Z_n|^2)^{1/2}\sim \frac{1}{4\sqrt{-g^{\prime}(0)}}. $$
Hence
\begin{align}\label{eqproofselectionbias}
& a_nn^{-\gamma/2}= -4g^{\prime}(0)\frac{1}{n}\sum_{m=1}^n\frac{m\ep|Z_m|-(m-1)\ep|Z_{m-1}|}{m^{\gamma/2}}\nonumber \\
\sim &-4g^{\prime}(0) \frac{1}{n}\left[\frac{n\ep|Z_n|}{n^{\gamma/2}}+\sum_{m=1}^{n-1} m\ep|Z_m|\left(\frac{1}{m^{\gamma/2}}-\frac{1}{(m+1)^{\gamma/2}}\right)\right]\\
\le & -4g^{\prime}(0)\frac{1}{n}\left[\frac{n(\ep|Z_n|^2)^{1/2}}{n^{\gamma/2}}+\sum_{m=1}^{n-1} m(\ep|Z_m|^2)^{1/2}\left(\frac{1}{m^{\gamma/2}}-\frac{1}{(m+1)^{\gamma/2}}\right)\right]\nonumber \\
\sim & \sqrt{|g^{\prime}(0)|}\left[n^{-\gamma/2}+\frac{1}{n}\sum_{m=1}^{n-1} m \left(\frac{1}{m^{\gamma/2}}-\frac{1}{(m+1)^{\gamma/2}}\right)\right]\sim \frac{2\sqrt{|g^{\prime}(0)|}}{2-\gamma} n^{-\gamma/2}.\nonumber
\end{align}
  Hence $\limsup\limits_{n\to \infty} a_n\le  2\sqrt{ |g^{\prime}(0)|}/ (2-\gamma)$.

  On the other hand,
 \begin{align*}
 \ep|Z_n|^3=   \frac{1}{n}\sum_{m=1}^n \sum_{\bm t} p(\bm t) \frac{\ep|\Lambda_m(\bm t)|^3}{m^{3\gamma/2}} =O(1)
  \end{align*}
and $\ep Z_n^2\le (\ep |Z_n|)^{1/2}(\ep |Z_n|^{3})^{1/2}$.
  It follows that
\begin{align*}
 \liminf\limits_{n\to \infty} a_n=& \liminf_{n\to\infty}\left[-4g^{\prime}(0)\right]\liminf\limits_{n\to \infty}\ep |Z_n|\\
  \ge & \liminf\limits_{n\to \infty}\left[-4g^{\prime}(0)\right]\left(\frac{\ep Z_n^2}{\ep Z_n^3} \right)^2\ge c_0 >0
  \end{align*}
by \eqref{eqproofselectionbias}. Thus
\begin{equation}\label{selectionBias}
SB_n =\frac{1}{2}+ O(n^{-\gamma/2}).
\end{equation}
is proved. To show the constant $c_0$ does not depend on $\gamma$. It is sufficient to show that for any integer $r\ge 0$, there is a constant $A_r$ which does not depend on $\gamma$ such that
\begin{equation}\label{uniformboundMoment}
\limsup_{n\to \infty} \frac{\ep M_n^r}{n^{\gamma r}}\le A_r,
\end{equation}
which implies that $\limsup_{n\to \infty}\ep[|Z_n|^3]\le C$ for a constant $C$ that does not depend on $\gamma$. 
We use  induction, assume \eqref{uniformboundMoment} holds for $r$. We now consider $r+1$. Let $m=m_n$ such that \eqref{eqinduction1} and \eqref{eqinduction2} hold. Similar to \eqref{TaylaExp2}, we have
\begin{align*}
\frac{c_r}{4(r+1)}& \ep(M_m+1)^r
\ge  \ep\left[(M_m+1)^rS_{m+1}(\bm\Lambda_m)\right]\\
=&-4g^{\prime}(0)\ep\left[(M_m+1)^r\frac{\sum_{\bm t}\Lambda_m^2(\bm t)p(\bm t)}{m^{\gamma}}\right]\\
&+ o(1)\left[\ep\left((M_m+1)^r\frac{\sum_{\bm t}\Lambda_m^2(\bm t)p(\bm t)}{m^{\gamma}}\right)^2\right]^{1/2}\\
\ge & c m^{-\gamma} \ep\left[(M_m+1)^rM_m\right]+ o(1) m^{-\gamma} \left[\ep M_m^{2(r+1)}\right]^{1/2}\\
=&c m^{-\gamma} \ep\left[M_m^{r+1}\right]+ o(1) m^{-\gamma} m^{\gamma(r+1)}.
\end{align*}
Hence, by \eqref{eqinduction1},
$$  \ep[M_n^{r+1}] \le m^{\gamma}\left(\frac{C_r}{4(r+1)c}+C_r\right)\ep(M_m+1)^r+o(1)   m^{\gamma(r+1)}. $$
It follows that
$$ \limsup_{n\to \infty} \frac{\ep[M_n^{r+1}]}{n^{\gamma(r+1)}}\le \left(\frac{C_r}{4(r+1)c}+C_r\right) A_r, $$
by the induction. This completes the proof.
   $\Box$
\smallskip

{\bf Proof of Corollary 3}. Note $ C^{-1} M_n\le \|\bm \Lambda_n\|^2\le C  M_n$. By the condition given in the corollary, when $4|\Lambda_n(\bm t)|\ge x_0$,
$$\left|g_n(4\Lambda_n(\bm t))-\frac{1}{2}\right|\ge p_0>0. $$
Let $E=\{4|\Lambda_m(\bm t)|<x_0\}$. Then on $E^c$,
$$S_{m+1}(\bm \Lambda_m)\ge c_0|\Lambda_m(\bm t)|. $$
It follows that
$$ \ep[|\Lambda_m(\bm t)|^{2r+2}E]\le  C, $$
and
$$\ep[|\Lambda_m(\bm t)|^{2r+1}E^c]\le  C     \ep\left[|\Lambda_m(\bm t)|^{2r}S_{m+1}(\bm\Lambda_m)\right] \le C \ep\left[(M_m+1)^rS_{m+1}(\bm\Lambda_m)\right]. $$
Let $m\in \{0,1,\ldots,n\}$ satisfy \eqref{eqinduction1} and \eqref{eqinduction2}. Then \eqref{eqinductionsad} remains true with $\gamma=0$. The remaining proof is similar to that of \eqref{eqth2.1} and \eqref{eqth2.2} by repeating  the induction  procedure \eqref{eqinduction3}-\eqref{eqinduction5} $i$ times until $1/2^{i-1}<\epsilon/2$. $\Box$

\section{Deriving the asymptotic selection bias}\label{sect:proof:SB}
\setcounter{equation}{0}

In this section we give the proofs of Theorem 3 in Section 3 and Remark 2 in Section 4.

Theorem 3 is a special case of Remark 2. 

   {\bf Proof of the Remark 2}. The proof of that $\sup_n\ep[M_n^r]<\infty$  can be found in Hu and Zhang (2020). For  completeness  of this paper, we give a direct  proof here. The allocation function probability $g(4\Lambda_{n-1}(\bm t))$ is now only a function of $\bm \Lambda_{n-1}$, so the Markov chain $\bm \Lambda_n$  is homogeneous. Note
   $0<g(x)<1$ for all $x$. The transition probability in  \eqref{tran:probab} is always positive. So the Markov chain is  irreducible.     It can be verified that
 $S_n(\bm \Lambda_{n-1})\to +\infty$ as $M_{n-1}\to \infty$.   By  \eqref{eq:conditionalM}, it follows that there are positive constants $b_r$ and $c_r$ such that
   $$
 \ep\big[M_n^{r+1}|\mathscr{F}_{n-1}\big]- M_{n-1}^{r+1}
\le   - \big(M_{n-1}+1\big)^r +b_r\mathbb I\left\{M_{n-1}\le c_r\right\}.
$$
Notice that $M_n=M^{\ast}(\bm\Lambda_n)$ is a    function of $\bm \Lambda_n$. The above inequality is just 
\begin{equation} \label{eq:drift} P_{\lambda}\big(M^{\ast}\big)^{r+1}-\big(M^{\ast}\big)^{r+1}\le -\big(M^{\ast}+1\big)^r+b_r\mathbb I\left\{M^{\ast}\le c_r\right\},
\end{equation}
where $M^{\ast}=M^{\ast}(\bm \Lambda)$, 
$(P_{\lambda} f)(\bm \Lambda) = \sum_{\bm \Lambda'\in\bm L(\mathbb{Z}^{m})}P_{\lambda}(\bm \Lambda,\bm \Lambda')f(\bm \Lambda') $
is the transition expectation of the function $f$,  $P_{\lambda}(\bm \Lambda,\bm \Lambda')$ is the transition probability from $\bm \Lambda$ to $\bm \Lambda'$. Then \eqref{eq:drift} implies that the irreducible Markov chain $\bm \Lambda_n$  with period 2 is positive recurrent with an invariant distribution $\bm \pi$ and $\sup_n\ep\big[(M_n+1)^r\big]<\infty$. In fact, write
$$ \overline{P}_{\lambda}^{(n)}(\bm \Lambda,\bm \Lambda')=\frac{1}{n}\sum_{m=0}^{n-1} P_{\lambda}^{(m)}(\bm \Lambda,\bm \Lambda'),$$
where $P_{\lambda}^{(m)}(\bm \Lambda,\bm \Lambda')$ is the $m$-step transition probability from    $\bm \Lambda$ to $\bm \Lambda'$. \eqref{eq:drift} implies that for any initial distribution $\bm \pi_0$,
\begin{align} \label{eq:drift2}
\limsup_{n\to \infty}  \bm \pi_0 \overline{P}_{\lambda}^{(n)}\big(M^{\ast}+1\big)^r\le & b_r \limsup_{n\to \infty}  \bm \pi_0\overline{P}_{\lambda}^{(n)}I\big\{M^{\ast}\le c_r\big\}+\limsup_{n\to \infty}\frac{\bm \pi_0\big(M^{\ast}\big)^{r+1}}{n} \nonumber \\
\le &  b_r \limsup_{n\to \infty}  \bm \pi_0 \overline{P}_{\lambda}^{(n)}I\big\{M^{\ast}\le c_r\big\} ,\; \forall r>0.
\end{align}
If the Markov-Chain is not positive recurrent, then for any $\bm \Lambda$ and $\bm \Lambda'$,  $ \overline{P}_{\lambda}^{(n)}(\bm \Lambda,\bm \Lambda')\to 0$. For each $C>0$, $\{M^{\ast}\le C\}$ contains only a finite number of values of the Markov chain. So,
$\bm \pi_0\overline{P}_{\lambda}^{(n)} I\left\{M^{\ast}\le C\right\}\to 0$.  It follows that the limit in \eqref{eq:drift2} is zero,
which is  obviously a contradiction since $\bm \pi_0 \overline{P}_{\lambda}^{(n)}\big(M^{\ast}+1\big)^r\ge 1$.  Hence $\bm\Lambda_n$ is a positive recurrent  Markov Chain and has an invariant distribution $\bm\pi$.  Now, \eqref{eq:drift2} implies
$$ \ep_{\pi}\left[\big(M^{\ast}+1\big)^r\right]=\bm \pi\overline{P}_{\lambda}^{(n)}\big(M^{\ast}+1\big)^r\le b_r, \; \forall r>0, $$
  which implies $\sup_n\ep[M_n^r]<\infty$ since $\frac{1}{2}\Big(\ep[M_{2n}^r]+\ep[M_{2n+1}^r]\Big)\to \ep_{\bm\pi}[(M^{\ast})^r]$.
  
Next, consider the selection bias $SB_n$. Notice that
\begin{align}\label{eq:proofofSB}
SB_n=&\frac{1}{n}\sum_{m=1}^n \ep\left[p_m\vee(1-p_m)\right]=\frac{1}{2}+ \frac{1}{n}\sum_{m=1}^n \ep\Big[\big|p_m-\frac{1}{2}\big|\Big]\nonumber\\
=&
\frac{1}{2}+ \frac{1}{n}\sum_{m=1}^n \sum_{\bm t}\ep\Big|g\big(4 \Lambda_{m-1}(\bm t) \big)-\frac{1}{2}\Big|p(\bm t)\nonumber\\
\to & \frac{1}{2}+\sum_{\bm t} \ep_{\bm \pi}\Big[\Big|g\big(4\Lambda(\bm t)\big)-\frac{1}{2}\Big|\Big]p(\bm t)>1/2,
\end{align}
 by the ergodic theorem. Here  the last inequality is due to $g(x)\not\equiv \frac{1}{2}$ and the fact that the Markov chain $\bm\Lambda_n$ is irreducible and so each possible state has a positive probability mass under the invariant distribution $\bm \pi$. 
This completes the proof. $\Box$.

\section{Further discussion on the randomness}\label{sect:proof:ran}
\setcounter{equation}{0}

{\em Selection bias under partial information guessing strategy}.  When the clinical trials  have  many covariates, one may think that it is hard for the experimenter  to know all the information so that he/or she can use an optimal guessing strategy, and so the selection bias can be ignored. We can show that this is not  true at least for the Pocock and Simon's (1975)  procedure. Write the allocation probability $p_m=g(\Lambda_{m-1}( \tilde{\bm X}_m))$ as a function of $  \Lambda_{m-1} ( \tilde{\bm X}_m)$. The guessing strategy is to guess treatment 1 when $G_m>0$, treatment 2 when $G_m<0$, and guess treatment 1 or 2 with probability $1/2$ when $G_m=0$, where $G_m$ is the guessing factor. Then  the expected proportion of correct guesses is
 \begin{align*}
 SB_n(g,G)=\frac{1}{2} +\frac{1}{n}\sum_{m=1}^n \ep\left[\big(p_m-\frac{1}{2}\big)\text{sgn}(G_m)\right],
 \end{align*}
 where "g" stands for the allocation $g(x)$, "G" stands for the guessing strategy $G_m$.  Because   the experimenter uses part or all information of $\Lambda_{m-1}(  \tilde{\bm X}_m )$, so  $G_m=G(\Lambda_{m-1}(  \tilde{\bm X}_m ))$ is a function of $\Lambda_{m-1}(  \tilde{\bm X}_m )$, where $G(\cdot)$ is a random function.
 We assume that under the guessing strategy, the allocation is  guessed in a correct  direction, i.e., for any point $x,y$,
 $\pr\big(\{G(x)-G(y)\}\cdot\{ g(x)- g(y)\}\ge 0\big)=P(x,y)\ge 1/2$, and there is at least a pair of points $x_0$ and $y_0$ such that $\pr\big(G(x_0)>0>G(y_0), g(x_0)> g(y_0)\big)>1/2$.
 \begin{theorem} For the Pocock and Simon's procedure, Hu and Hu's procedure or more general procedures as in Remark 2, we have
 $$ \lim_{n\to \infty} SB_n(g,G)>\frac{1}{2}. $$
 \end{theorem}
 {\bf Proof.} Notice  Remark 2. The sequence $\{\Lambda_n(\bm t); t_j=1,\cdots,m_j, j=1,\cdots, I\}_{n=1}^{\infty}$ is  an    irreducible   positive recurrent multi-dimensional Markov Chain having an invariant  distribution $\bm \pi$.  Let $\Lambda^{\ast}(\bm t), \tilde{\bm X}^{\ast}$ be independent copies of $\Lambda(\bm t), \tilde{\bm X}$ and define  $\eta=\Lambda( \tilde{\bm X})$,  $\eta^{\ast}=\Lambda^{\ast}( \tilde{\bm X}^{\ast})$.
 Then  by the ergodic theorem,
  \begin{align*}
  &SB_n(g,G)   -\frac{1}{2}= \frac{1}{n}\sum_{m=1}^n \ep\left[\big(p_m-\frac{1}{2}\big)\text{sgn}(G_m)\right] \\
  \longrightarrow  & \ep_{\bm \pi}\Big[ \Big(g\big(\eta)\big)-\frac{1}{2}\Big) \text{sgn}[G(\eta)]\Big]\\
  =   &\ep_{\bm \pi}\Big[ g\big(\eta)\big)-\frac{1}{2} \Big]\ep\left[\text{sgn}[G(\eta)]\right]+\\
  &\frac{1}{2}\ep_{\bm \pi\times \bm \pi}\left[ \big(g(\eta)-g(\eta^{\ast})\big) \Big(\text{sgn}(G(\eta))-\text{sgn}(G(\eta^{\ast}))\Big) \right]\\
   =  & \ep_{\bm \pi}\Big[ g\big(\eta)\big)-\frac{1}{2} \Big]\ep\left[\text{sgn}[G(\eta)]\right]\\
 &  +\frac{1}{2}\ep_{\bm \pi\times \bm \pi}\left[\left|\big(g(\eta)-g(\eta^{\ast})\big) \Big(\text{sgn}(G(\eta))-\text{sgn}(G(\eta^{\ast}))\Big)\right|
 \left(2P(\eta,\eta^{\ast})-1\right)\right].
 \end{align*}
 For the value of the limit, the second   term is positive due to the assumption for guessing in a correct  direction. The first term is zero because $\ep_{\bm \pi}g(\eta)=1/2$ under the invariant distribution $\bm\pi$, which can be verified by
$$0\leftarrow\frac{\ep_{\bm\pi}[D_n]}{n}=\frac{1}{n}\sum_{j=1}^n \ep_{\bm\pi}[2T_j-1]=\frac{1}{n}\sum_{j=1}^n 2(\ep_{\bm\pi} p_j-1)=2\ep_{\bm\pi}g(\eta)-1. $$
  This completes the proof. $\Box$

{\em Entropy of a CAR procedure}. 
Besides the selection bias,  another measure of the randomness is the entropy.  The entropy of a probability $p$ is defined by
$$ E( p)=-  p\log p-(1-p)\log (1-p). $$
The entropy of the design with allocation probabilities $p_{m}$s is  defined by
$$Ent_n=\ep\left[\frac{1}{n}\sum_{m=1}^n E(p_{m})\right]. $$
It is obvious that $0\le Ent_n\le \log 2$. $Ent_n$ takes its largest value $\log 2$ when and only when the randomization is  complete randomization in which $p_{m}=1/2$ for all  $m$. $Ent_n$ takes its smallest value $0$ when and only when each patient is assigned deterministically.

\begin{theorem} We have
$$ Ent_n \to \log 2\; \text{ if and only if }\; SB_n\to \frac{1}{2}. $$
\end{theorem}

{\bf Proof.} Since $SB_n=\frac{1}{2}+ \frac{1}{n}\sum_{m=1}^n\ep\big|p_m-\frac{1}{2}\big|$, it is sufficient to notice that 
$$ 2\Big(p-\frac{1}{2}\Big)^2\le E\Big(\frac{1}{2}\Big)-E(p) \le (2\log 2)\Big|p-\frac{1}{2}\Big|,$$
and
$$ \frac{1}{n}\sum_{m=1}^n\ep\big|p_m-\frac{1}{2}\big|
\le \left(\frac{1}{n}\sum_{m=1}^n\ep \big|p_m-\frac{1}{2}\big|^2\right)^{1/2}. \qquad \Box$$

 \begin{theorem} For the Pocock and Simon's procedure, Hu and Hu's procedure or more general procedures as in Remark 2, we have
 $$ \lim_{n\to \infty} Ent_n<\log 2. $$
 \end{theorem}
 {\bf Proof}. Similar to \eqref{eq:proofofSB}, we have
\begin{align*}
Ent_n=&\frac{1}{n}\sum_{m=1}^n \ep\left[E(p_m)\right]
=
  \frac{1}{n}\sum_{m=1}^n \sum_{\bm t}\ep\left[ E\Big(g\big(4 \Lambda_{m-1}(\bm t)\big) \Big)\right]p(\bm t) \\
\to &  \sum_{\bm t} \ep_{\bm \pi}\left[ E\Big(g\big(4 \Lambda(\bm t) \big)\Big)\right]p(\bm t)<\log 2,
\end{align*}
 by the ergodic theorem. Here  the last inequality is due to
  $E(x)<\log 2$ when $x\ne 1/2$, $g(x)\not\equiv 1/2$ and the fact that   each possible state of $\bm\Lambda_n$ has a positive probability mass under the invariant distribution $\bm \pi$. $\Box$

\section{Simulation Study}\label{sect:simu}

\setcounter{equation}{0}

In this section, we study the new CAR by simulation. 
\begin{description}
	\item[\sl Model Setting:] The underlying model is a linear model with independent covariates:
	$$ Y_i=T_i\mu_1+(1-T_i)\mu_2+\sum_{j=1}^3 \beta_j X_{i,j}+\epsilon_i,$$
	where $\beta_j=1$ for $j = 1, 2,3$, $X_{i,1}\sim N(0, 1)$, $X_{i,2}, X_{i,3}\sim N(1, 1)$, and the covariates are mutually independent. The random error $\epsilon_i\sim N(0, 2^2)$ is independent of all $X_{i,j}$.
\end{description}

For randomization, we employ complete randomization (CR) and four covariate-adaptive randomization (CAR) procedures to assign patients to treatments: Pocock and Simon's procedure (PS) and the new CAR procedures (NEW) with $\gamma=0.2$, $0.5$, and $0.8$, based on the discretized covariates $\widetilde{\bm X}_i=(d(X_{i,1}), d(X_{i,2}), d(X_{i,3}))$, where $d(X)=0I\{X\le 0\}+1I\{0< X< 2\}+2I\{X\ge 2\}$. 
Pocock and Simon's procedure is the general CAR with imbalance measures
$$
Imb_{n}^{(t)}= \frac{1}{I}\sum_{i=1}^{I}[D^{(t)}_n (i;t_i^{\ast})]^2,
$$
and the allocation probability
\begin{align}\label{func_PS}
 p_n=g_n(Imb_{n}^{(1)}-Imb_{n}^{(2)})=\begin{cases}
	1-\rho, & Imb_{n}^{(1)}-Imb_{n}^{(2)} >0 , \\
	\frac{1}{2}, &  Imb_{n}^{(1)}-Imb_{n}^{(2)} = 0,\\
	\rho, &  Imb_{n}^{(1)}-Imb_{n}^{(2)} < 0,
\end{cases}
\end{align} 
where $1/2<\rho<1$. For the new CAR, we choose imbalance measures the same as that in Pocock and Simon's procedure, and the function $g$ in (4.8) to be $g(x)=(1-\rho)\vee\big(\frac{1}{2}-\frac{1}{2} \operatorname*{sgn}(x)(1\wedge |x|)\big)\wedge \rho$, $0< \gamma< 1$ and $1/2<\rho<1$, i.e.,  the allocation probability of the $n$th patients is 
\begin{align}\label{func_new}
&p_n=g_n(Imb_{n}^{(1)}-Imb_{n}^{(2)}) \nonumber \\
=&(1-\rho)\vee\left[\frac{1}{2}-\frac{1}{2}\operatorname*{sgn}(Imb_{n}^{(1)}-Imb_{n}^{(2)}) \Big(1\wedge \frac{|Imb_{n}^{(1)}-Imb_{n}^{(2)}|}{(n-1)^{\gamma}}\Big)  \right]\wedge \rho.
\end{align} 
This CAR procedure has the property that, when $|Imb_{n}^{(1)}-Imb_{n}^{(2)}|/(n-1)^{\gamma}$ is large, the allocation probability is the same as that of Pocock and Simon's procedure, and when 
$|Imb_{n}^{(1)}-Imb_{n}^{(2)}|/(n-1)^{\gamma}$ is small, the allocation probability is a linear function of $(Imb_{n}^{(1)}-Imb_{n}^{(2)})/(n-1)^{\gamma}$:
$$p_n=g_n(Imb_{n}^{(1)}-Imb_{n}^{(2)})=\frac{1}{2}-\frac{1}{2} \frac{Imb_{n}^{(1)}-Imb_{n}^{(2)}}{(n-1)^{\gamma}}. $$

Under each randomization procedure, we apply the two-sided \(t\)-test for
$$ H_0: \mu_1=\mu_2 \quad \text{vs.} \quad H_A: \mu_1\ne \mu_2, $$
using the test statistic defined in Equation (2.3) of the main text.

We evaluate the Type I error rates under $H_0$, power under $H_A:\mu_1-\mu_2=\delta/\sqrt{n}$ for $\delta=5,10,15$, and selection bias for sample sizes $n\in\{50,100,200,400\}$. Each simulation result is based on 5,000 replications. 

The results are reported in Table \ref{Simulation}. For complete randomization, the selection bias is fixed at 0.5, obtained analytically rather than through simulation.

\begin{table}[htbp]
	\centering
	\footnotesize 
	\renewcommand{\arraystretch}{0.7}
	\caption{Simulation results for selection bias, Type I error rates, and power under various randomization procedures and sample sizes.}
	\begin{tabular}{*{8}{c}}
		\toprule
		$n$ & Ran &Allocation Probability & SB & Type I & \multicolumn{3}{c}{Power}\\
		\cmidrule(lr{0pt}){6-8} 
		& & & & & $\delta=5$ & $\delta=10$ & $\delta=15$ \\
		\midrule
	50	 & CR & - 
		& 0.5000 & 0.0530 & 0.2120 & 0.6440 & 0.9356 \\
\cmidrule(lr{0pt}){2-8}
		& PS & \eqref{func_PS} with $\rho=0.9$ 
		& 0.8334 & 0.0550 & 0.2260 & 0.6672 & 0.9500 \\
		& NEW & \eqref{func_new} with $\gamma=0.2,\rho=0.9$ 
		& 0.8186 & 0.0486 & 0.2328 & 0.6686 & 0.9502 \\
		& NEW & \eqref{func_new} with $\gamma=0.5,\rho=0.9$ 
		& 0.7391 & 0.0500 & 0.2280 & 0.6648 & 0.9496 \\
		& NEW & \eqref{func_new} with $\gamma=0.8,\rho=0.9$ 
		& 0.6544 & 0.0510 & 0.2284 & 0.6696 & 0.9466 \\
		\cmidrule(lr{0pt}){2-8} 
		& PS & \eqref{func_PS} with $\rho=2/3$ 
		& 0.6444 & 0.0496 & 0.2222 & 0.6646 & 0.9468 \\
		& NEW & \eqref{func_new} with $\gamma=0.2,\rho=2/3$ 
		& 0.6444 & 0.0506 & 0.2254 & 0.6664 & 0.9442 \\
		& NEW & \eqref{func_new} with $\gamma=0.5,\rho=2/3$ 
		& 0.6386 & 0.0518 & 0.2312 & 0.6596 & 0.9476 \\
		& NEW & \eqref{func_new} with $\gamma=0.8,\rho=2/3$ 
		& 0.6156 & 0.0510 & 0.2148 & 0.6696 & 0.9486 \\
		\midrule

	100	& CR & - 
		& 0.5000 & 0.0620 & 0.2328 & 0.6776 & 0.9490 \\
\cmidrule(lr{0pt}){2-8}
		 & PS & \eqref{func_PS} with $\rho=0.9$ 
		& 0.8254 & 0.0482 & 0.2330 & 0.6878 & 0.9578 \\
		& NEW & \eqref{func_new} with $\gamma=0.2,\rho=0.9$ 
		& 0.8051 & 0.0506 & 0.2348 & 0.6950 & 0.9580 \\
		& NEW & \eqref{func_new} with $\gamma=0.5,\rho=0.9$ 
		& 0.7098 & 0.0504 & 0.2304 & 0.7032 & 0.9562 \\
		& NEW & \eqref{func_new} with $\gamma=0.8,\rho=0.9$ 
		& 0.6209 & 0.0468 & 0.2356 & 0.6882 & 0.9576 \\
		\cmidrule(lr{0pt}){2-8} 
		& PS & \eqref{func_PS} with $\rho=2/3$ 
		& 0.6493 & 0.0504 & 0.2380 & 0.6966 & 0.9568 \\
		& NEW & \eqref{func_new} with $\gamma=0.2,\rho=2/3$ 
		& 0.6491 & 0.0512 & 0.2354 & 0.6922 & 0.9534 \\
		& NEW & \eqref{func_new} with $\gamma=0.5,\rho=2/3$ 
		& 0.6384 & 0.0486 & 0.2304 & 0.6916 & 0.9598 \\
		& NEW & \eqref{func_new} with $\gamma=0.8,\rho=2/3$ 
		& 0.6026 & 0.0508 & 0.2382 & 0.6894 & 0.9558 \\
		\midrule

	200	& CR & - 
		& 0.5000 & 0.0492 & 0.2366 & 0.6910 & 0.9554 \\
\cmidrule(lr{0pt}){2-8}
		& PS & \eqref{func_PS} with $\rho=0.9$ 
		& 0.8336 & 0.0462 & 0.2424 & 0.7024 & 0.9616 \\
		& NEW & \eqref{func_new} with $\gamma=0.2,\rho=0.9$ 
		& 0.7994 & 0.0452 & 0.2298 & 0.7160 & 0.9626 \\
		& NEW & \eqref{func_new} with $\gamma=0.5,\rho=0.9$ 
		& 0.6856 & 0.0506 & 0.2382 & 0.6986 & 0.9608 \\
		& NEW & \eqref{func_new} with $\gamma=0.8,\rho=0.9$ 
		& 0.5960 & 0.0536 & 0.2426 & 0.7012 & 0.9608 \\
		\cmidrule(lr{0pt}){2-8} 
		& PS & \eqref{func_PS} with $\rho=2/3$ 
		& 0.6512 & 0.0520 & 0.2286 & 0.6880 & 0.9590 \\
		& NEW & \eqref{func_new} with $\gamma=0.2,\rho=2/3$ 
		& 0.6506 & 0.0472 & 0.2334 & 0.7004 & 0.9590 \\
		& NEW & \eqref{func_new} with $\gamma=0.5,\rho=2/3$ 
		& 0.6332 & 0.0508 & 0.2310 & 0.6938 & 0.9570 \\
		& NEW & \eqref{func_new} with $\gamma=0.8,\rho=2/3$ 
		& 0.5866 & 0.0520 & 0.2320 & 0.6908 & 0.9570 \\
		\midrule

	400	& CR & - 
		& 0.5000 & 0.0498 & 0.2308 & 0.6914 & 0.9590 \\
\cmidrule(lr{0pt}){2-8}
		& PS & \eqref{func_PS} with $\rho=0.9$ 
		& 0.8295 & 0.0576 & 0.2376 & 0.7082 & 0.9612 \\
		& NEW & \eqref{func_new} with $\gamma=0.2,\rho=0.9$ 
		& 0.7903 & 0.0524 & 0.2408 & 0.7026 & 0.9614 \\
		& NEW & \eqref{func_new} with $\gamma=0.5,\rho=0.9$ 
		& 0.6596 & 0.0492 & 0.2384 & 0.7002 & 0.9628 \\
		& NEW & \eqref{func_new} with $\gamma=0.8,\rho=0.9$ 
		& 0.5743 & 0.0516 & 0.2356 & 0.6994 & 0.9588 \\
		\cmidrule(lr{0pt}){2-8} 
		& PS & \eqref{func_PS} with $\rho=2/3$ 
		& 0.6531 & 0.0504 & 0.2260 & 0.6918 & 0.9618 \\
		& NEW & \eqref{func_new} with $\gamma=0.2,\rho=2/3$ 
		& 0.6526 & 0.0462 & 0.2298 & 0.7056 & 0.9642 \\
		& NEW & \eqref{func_new} with $\gamma=0.5,\rho=2/3$ 
		& 0.6256 & 0.0502 & 0.2398 & 0.7116 & 0.9616 \\
		& NEW & \eqref{func_new} with $\gamma=0.8,\rho=2/3$ 
		& 0.5694 & 0.0528 & 0.2400 & 0.7056 & 0.9604 \\
		\bottomrule
	\end{tabular}
    \label{Simulation}
\end{table}

First, across all settings, the Type I error rates are close to the nominal level of 5\%, indicating that the hypothesis tests are well calibrated. Second, for relatively small sample sizes, both the new CAR procedures and the PS procedure yield slightly higher power than complete randomization; however, this power advantage gradually diminishes as the sample size increases. Third, for the same value of \(\rho\), the selection bias under the PS procedure is consistently larger than that under the new CAR procedures. Furthermore, smaller values of \(\rho\) lead to lower selection bias for both the new CAR and PS procedures, and the larger values of $\gamma$ lead to lower selection bias for the new CAR. Notably, the selection bias of the PS procedure does not decrease with increasing sample size, whereas the new CAR procedures exhibit a decreasing selection bias as the sample size grows. In particular, the reduction is more pronounced for larger values of \(\gamma\), indicating a faster convergence rate. This observation aligns with our theoretical results.

\section*{Acknowledgments}
We are grateful to Mr. Yuhang Tao, a Ph.D. student under the author's supervision, for contributing the simulation data in Table \ref{Simulation}.

\end{document}